\documentclass[12pt,leqno]{article}
\usepackage{amsfonts}
\usepackage{amsmath, amsthm, amsfonts, amssymb, color}
\usepackage{mathrsfs}
\usepackage{color}
\theoremstyle{definition}

\newcommand{\scr}[1]{\mathscr #1}
\definecolor{wco}{rgb}{0.5,0.2,0.3}

\numberwithin{equation}{section} \theoremstyle{remark}

\newcommand{\ua}{\uparrow}

\title{{\bf Distribution Dependent SDEs for Landau Type Equations}\footnote{Supported in
 part by  NNSFC (11431014,11626245,11626250).} }
\author{
{\bf     Feng-Yu Wang  }\\
\footnotesize{Center for Applied Mathematics, Tianjin University, Tianjin 300072, China}\\
 \footnotesize{ Department of Mathematics,
Swansea University, Singleton Park, SA2 8PP, United Kingdom}\\
\footnotesize{  wangfy@tju.edu.cn, F.-Y.Wang@swansea.ac.uk}}
\begin{document}
\allowdisplaybreaks
\def\R{\mathbb R}  \def\ff{\frac} \def\ss{\sqrt} \def\B{\mathbf
B} \def\W{\mathbb W}
\def\N{\mathbb N} \def\kk{\kappa} \def\m{{\bf m}}
\def\ee{\varepsilon}\def\ddd{D^*}
\def\dd{\delta} \def\DD{\Delta} \def\vv{\varepsilon} \def\rr{\rho}
\def\<{\langle} \def\>{\rangle} \def\GG{\Gamma} \def\gg{\gamma}
  \def\nn{\nabla} \def\pp{\partial} \def\E{\mathbb E}
\def\d{\text{\rm{d}}} \def\bb{\beta} \def\aa{\alpha} \def\D{\scr D}
  \def\si{\sigma} \def\ess{\text{\rm{ess}}}
\def\beg{\begin} \def\beq{\begin{equation}}  \def\F{\scr F}
\def\Ric{\text{\rm{Ric}}} \def\Hess{\text{\rm{Hess}}}
\def\e{\text{\rm{e}}} \def\ua{\underline a} \def\OO{\Omega}  \def\oo{\omega}
 \def\tt{\tilde} \def\Ric{\text{\rm{Ric}}}
\def\cut{\text{\rm{cut}}} \def\P{\mathbb P} \def\ifn{I_n(f^{\bigotimes n})}
\def\C{\scr C}      \def\aaa{\mathbf{r}}     \def\r{r}
\def\gap{\text{\rm{gap}}} \def\prr{\pi_{{\bf m},\varrho}}  \def\r{\mathbf r}
\def\Z{\mathbb Z} \def\vrr{\varrho} \def\ll{\lambda}
\def\L{\scr L}\def\Tt{\tt} \def\TT{\tt}\def\II{\mathbb I}
\def\i{{\rm in}}\def\Sect{{\rm Sect}}  \def\H{\mathbb H}
\def\M{\scr M}\def\Q{\mathbb Q} \def\texto{\text{o}} \def\LL{\Lambda}
\def\Rank{{\rm Rank}} \def\B{\scr B} \def\i{{\rm i}} \def\HR{\hat{\R}^d}
\def\to{\rightarrow}\def\l{\ell}\def\iint{\int}
\def\EE{\scr E}\def\Cut{{\rm Cut}}
\def\A{\scr A} \def\Lip{{\rm Lip}}
\def\BB{\scr B}\def\Ent{{\rm Ent}}\def\L{\scr L}

\maketitle

\begin{abstract} The distribution dependent stochastic differential equations  (DDSDEs) describe stochastic systems whose evolution is   determined by both   the  microcosmic site  and the  macrocosmic  distribution of the particle. The density function associated with a DDSDE solves a  nonlinear PDE.    Due to the distribution dependence,  some standard techniques developed for  SDEs do not apply.  By iterating in distributions,  a strong solution is constructed using SDEs with control. By proving the uniqueness, the distribution of solutions is identified with a nonlinear semigroup $P_t^*$ on the space of probability measures. The  exponential   contraction   as well as  Harnack inequalities and applications  are investigated for the nonlinear semigroup $P_t^*$   using coupling by change of measures. The main results are illustrated by  homogeneous  Landau  equations.
\end{abstract} \noindent
 \vskip 2cm

\section{Introduction}

A fundamental application of the It\^o SDE is to solve  Kolmogorov's problem \cite{Kol} of determining Markov processes whose distribution density satisfies the Fokker-Planck-Kolmogorov equation. Let $W_t$ be the $d$-dimensional Brownian motion on a complete  probability space with natural filtration $(\OO, \{\F_t\}_{t\ge 0}, \P)$, and let $b:  \R^d\to \R^d; \si: \R^d\to \R^d\otimes \R^d$ be measurable. Then   the distribution density of the solution to the SDE
\beq\label{E0} \d X_t= b(X_t)\d t+ \si(X_t)\d W_t\end{equation} satisfies the parabolic equation
\beq\label{LH}\pp_t f_t = \ff 1 2 \sum_{i,j=1}^d \pp_i\pp_j\big\{(\si\si^*)_{ij}  f_t \big\} -\sum_{i=1}^d \pp_i\{b_if_t\},\end{equation}
which describes the time evolution of the probability density function of the velocity of a particle under the influence of drag forces and random forces.
 If $b$ and $\si$ are $``$almost" locally Lipschitzian,   then the SDE \eqref{E0} has a unique strong solution up to life time
(c.f. \cite{FZ}). When $\si$ is invertible (i.e. the SDE is non-degenerate),   this condition has been largely weakened as $|b|+|\nn \si|\in L_{loc}^p (\d x)$ for some $p>d$,  see \cite{Zhang} and references within.

When coefficients $\si$ and $b$ also depend on the distribution of the solution, the SDE is called distribution dependent. This type SDEs often arise from mathematical physics, see for instance \cite{FM} for distribution dependent SDEs of jump type describing the Boltzmann equation, and \cite{Carlen} for those of diffusion type in Nelson's stochastic mechanics. Consider,   for instance, the Landau type equation
\beq\label{LD} \pp_t f_t=\ff 1 2 \text{div}\bigg\{\int_{\R^d} a(\cdot-z)\big(f_t(z)\nn f_t-f_t\nn f_t(z)\big)\d z\bigg\},\end{equation}
for some reference coefficient $a:\R^d\to \R^d\otimes \R^d$. This includes  the homogenous Landau equation where $d=3$ and
\beq\label{LLD} a(x)= |x|^{2+\gg} \Big(I-\ff{x\otimes x}{|x|^2}\Big),\ \ x\in \R^3 \end{equation} for some constant $\gg\in [-3,1].$ Landau equation is $``$grazing collision limit" of the Boltzmann equation. 
 When $\gg \in [0,1]$, the existence, uniqueness, regularity estimates, and exponential convergence   have been investigated in \cite{DV1,DV2, CA} and references within for  initial density in $L^1_{s_1}\cap L^2_{s_2}$ for large enough $s_1,s_2>0$, where $f\in L_s^p$ means $\int_{\R^3} |f(x)|^p(1+|x|^2)^{\ff s 2} \d x<\infty.$ 
   To describe the solution of \eqref{LD} using stochastic processes,  consider
the following distribution dependent SDE (DDSDE) for   $b= \text{div} a$ and $ \si$ such that $\si\si^*= a$:
\beq\label{SD0} \d X_t= (b*\L_{X_t})(X_t)\d t +(\si *\L_{X_t})(X_t)\d W_t,\end{equation} where $\L_\xi$ denotes the distribution of a random variable $\xi$, and
$$(f*\mu)(x):=\int_{\R^d} f(x-z)\mu(\d z)$$ for a function $f$ and a probability measure $\mu$.  By It\^o's formula and the integration by parts formula, the distribution density of $X_t$ is a weak solution to \eqref{LD}. For the homogenous Landau equation  with $\gg\in [0,1]$ and initial distribution density $f_0$ satisfying
\beq\label{DST} \int_{\R^3} f_0(x)\big(f_0(x)+ \e^{|x|^\aa}\big)\d x<\infty\ \ \text{for\ some\ }\aa>\gg, \end{equation}  the existence and uniqueness of weak solutions to \eqref{SD0} has been proved  in \cite{FG} by an approximation argument using particle systems. This approximation is known as \emph{propagation of chaos} according to Kac \cite{Kac}, see also \cite{Fu, Gu} and references within.

In this paper, we aim to investigate the (pathwise) strong solutions of  \eqref{SD0} and characterize their distribution properties.

In general, for measurable maps
$$b: [0,\infty)\times \R^d\times \scr P(\R^d)\to \R^d;\ \ \si: [0,\infty)\times \R^d\times \scr P(\R^d)\to \R^d\otimes\R^d,$$
we consider the following DDSDE on $\R^d$:
\beq\label{E1} \d X_t= b_t(X_t,\L_{X_t})\,\d t+ \si_t(X_t, \L_{X_t})\,\d W_t.\end{equation}
When more than one probability measures on $\OO$ are concerned, we use $\L_{X_t}|_\P$ instead of $\L_{X_t}$ to emphasize the distribution under probability $\P$.
 Due to technical reasons, we will restrict ourselves to the following   subspace of $\scr P$  for some $\theta\in [1,\infty)$:
$$\scr P_\theta:=\bigg\{\nu\in \scr P:\ \nu(|\cdot|^\theta):=\int_{\R^d} |x|^\theta\nu(\d x)<\infty\bigg\},$$ which is a polish space under the
 $L^\theta$-Wasserstein distance
 $$\W_\theta(\mu_1,\mu_2):= \inf_{\pi\in \scr C(\mu_1,\mu_2)}\bigg(\int_{\R^d\times\R^d} |x-y|^\theta\pi(\d x,\d y)\bigg)^{\ff 1 \theta},\ \ \mu_1,\mu_2\in \scr P_\theta,$$
 where $\scr C(\mu_1,\mu_2)$ is the set of all couplings for $\mu_1$ and $\mu_2$.

 The following definition is standard in the literature of SDEs.

\beg{defn}  $(1)$ For any $s\ge 0$, a continuous adapted process $(X_t)_{t\ge s}$ on $\R^d$ is called a (strong) solution of \eqref{E1} from time $s$, if
  $$\int_s^t \E\big\{|b_r(X_r,\L_{X_r})|+\|\si_r(X_r, \L_{X_r})\|^2\big\}\d r<\infty,\ \ t>s,$$   and $\P$-a.s.,
$$X_t = X_s +\int_s^t b_r(X_r, \L_{X_r})\d r + \int_s^t \si_r(X_r,\L_{X_r})\d W_r,\ \ t\ge s.$$
We say that \eqref{E1} has (strong or pathwise) existence and uniqueness in $\scr P_\theta$, if for any $s\ge 0$ and $\F_s$-measurable random variable $X_{s,s}$ with $\E|X_{s,s}|^\theta<\infty$, the equation from time $s$ has a unique solution  $(X_{s,t})_{t\ge s}$ with $\E |X_{s,t}|^\theta<\infty$. We   simply denote $X_{0,t}=X_t$.

$(2)$ A couple $(\tt X_t, \tt W_t)_{t\ge s}$ is called a weak solution to \eqref{E1} from time $s$, if $\tt W_t$  is the $d$-dimensional Brownian motion  with respect to a complete filtration  probability space  $ (\tt\OO, \{\tt\F_t\}_{t\ge 0}, \tt\P)$, and   $\tt X_t$ solves the DDSDE
\beq\label{E1'} \d \tt X_t = b_t(\tt X_t,  \L_{\tt X_t}|_{\tt\P})\d t + \si_t(\tt X_t,  \L_{\tt X_t}|_{\tt\P})\d \tt W_t,\ \ t\ge s.\end{equation}

$(3)$  \eqref{E1} is said to have weak uniqueness in $\scr P_\theta$, if for any $s\ge 0$, any two weak solutions of the equation  from time $s$ with common initial  distribution in $\scr P_\theta$  are equal in law. Precisely,  if $s\ge 0$ and $(\bar X_{s,t},\bar W_t)_{t\ge s}$
with respect to $(\bar\OO,  \{\bar\F_t\}_{t\ge 0}, \bar\P)$ and $(\tt X_{s,t}, \tt W_t)_{t\ge s}$ with respect to $ (\tt\OO, \{\tt\F_t\}_{t\ge 0}, \tt\P)$ are weak solutions of \eqref{E1}, then $\L_{\bar X_{s,s}}|_{\bar\P}= \L_{\tt X_{s,s}}|_{\tt\P}$ implies $\L_{\bar X_{s,\cdot}}|_{\bar\P}= \L_{\tt X_{s,\cdot}}|_{\tt\P}$.
\end{defn}

 When \eqref{E1} has strong existence and uniqueness in $\scr P_\theta$,   the solution $(X_t)_{t\ge 0}$ is a Markov process in the sense  that for any $s\ge 0$,  $(X_t)_{t\ge s}$ is determined by solving the equation from time $s$ with initial state   $X_s$. More precisely, letting $\{X_{s,t}(\xi)\}_{t\ge s}$ denote the solution of the equation from time $s$ with initial state  $X_{s,s}=\xi$, the existence and uniqueness imply
 \beq\label{MK} X_{s,t}(\xi)= X_{r,t}(X_{s,r}(\xi)),\ \ t\ge r\ge  s\ge 0, \xi\ {\rm is}\ \F_s\text{-measurable\ with\ } \E|\xi|^\theta<\infty.\end{equation}
 However, in general the solution is  not strong Markovian since  we do not have $\L_{X_\tau}= \L_{X_t}$ on the set $\{\tau=t\}$ for a stopping time $\tau$ and $t>0$. Moreover,   the associated Markov operators $(P_{s,t})_{t\ge s}$ given by
 $$P_{s,t} f(x):= \E f(X_{s,t}(x)),\ \ x\in \R^d, f\in \B_b(\R^d)$$ is not a semigroup, see \eqref{NM} below.

When the DDSDE  has   $\scr P_\theta$-weak uniqueness  (in the classical case  it follows from the pathwise uniqueness according to Yamada-Watanabe),   we may define a semigroup $(P_{s,t}^*)_{t\ge s}$ on $\scr P_\theta$ by letting $P_{s,t}^*\mu=\L_{X_{s,t}}$ for   $\L_{X_{s,s}}=\mu\in \scr P_\theta$.  Indeed, by \eqref{MK} we have
\beq\label{SM} P_{s,t}^*= P_{r,t}^* P_{s,r}^*,\ \ t\ge r\ge  s\ge 0.\end{equation}
To see that $(P_{s,t})_{t\ge s}$ is not a semigroup,  we write
$$(P_{s,t} f)(\mu)= (P_{s,t}^*\mu)(f):=\int_{\R^d} f\d (P_{s,t}^*\mu),\ \  f\in \B_b(\R^d), t\ge 0, \mu\in \scr P_\theta.$$  Then $P_{s,t}f(x)= (P_{s,t}f)(\dd_x)$, where $\dd_x$ is the Dirac measure at point $x$.
Since $(\L_{X_{s,t}})_{t\ge s}$ solves a nonlinear equation as indicated in the beginning,   the semigroup $P_{s,t}^*$ is   nonlinear; i.e.
$$ P_{s,t}^*\mu \ne \int_{\R^d} (P_{s,t}^*\dd_x )\mu(\d x),\ \ t>s\ge 0  $$ for a non-trivial distribution $\mu$.
In other words,   in general
$$(P_{s,t} f)(\mu)\ne\mu(P_{s,t} f):=\int_{\R^d} P_{s,t} f \d\mu,\ \ t>s\ge 0,$$
so that
\beq\label{NM} \beg{split} (P_{s,t}f)(\mu)&:= \int_{\R^d} f\d (P_{s,t}^*\mu)= \int_{\R^d} f\d (P_{r,t}^*P_{s,r}^*\mu)=(P_{r,t}f)(P_{s,r}^*\mu)\\
&\ne \int_{\R^d} (P_{r,t}f)\d (P_{s,r}^*\mu)= (P_{s,r}(P_{r,t} f))(\mu),\ \ t>s\ge 0.\end{split}\end{equation}

Although the semigroup $P_{s,t}^*$ is nonlinear, we may also investigate the ergodicity in the time homogeneous case when $\si_t$ and $b_t$ do not depend on $t$. In this case we have $P_{s,t}^*= P_{t-s}^*$ for $t\ge s\ge 0$.   We call $\mu\in \scr P_\theta$ an invariant probability measure of $P_t^*$ if $ P_t^*\mu=\mu$ for all $t\ge 0$, and we call the solution ergodic if there exists $\mu\in \scr P_\theta$ such that $\lim_{t\to\infty}   P_t^*\nu = \mu$ weakly  for any $\nu\in \scr P_\theta$.   Obviously, the ergodicity implies that $P_t^*$ has a  unique   invariant probability measure.

When $b$ and $\si$ are bounded and Lipschitz continuous in $(x,\nu)\in \R^d\times \W_2,$ the weak solution of \eqref{E1} has been constructed in \cite{SV} by using propagation of chaos. In this paper, we investigate the existence, uniqueness and distribution properties of the strong solutions. To explain the  difficulty of the study, let us recall some standard techniques developed  for \eqref{E0}   with locally bounded coefficients. Firstly, by a truncation argument one reduces an SDE with locally bounded coefficients to that with bounded coefficients, so that  when $\si$ is invertible   the existence of weak solutions is ensured by the Girsanov transform and the uniqueness follows from Zvonkin type argument, see e.g. \cite{Zhang} and references within. Then   the SDE has a unique strong solution according to   Yamada-Watanabe's principle \cite{YW}. However, these techniques do not apply to DDSDEs:  since the coefficients depend on the distribution which is not pathwisely determined, the truncation argument  and  Yamada-Watanabe's principle do not work; since the distribution of solution depends on the reference probability measure, the Girsanov transform method is   invalid for the construction of weak solutions. Moreover,   due to the lack of strong Markovian property, one can not let the marginal processes move together after the coupling time, so that the classical coupling argument does not apply.  To overcome the difficulty caused by distribution dependence, we will approximate the DDSDE \eqref{E1} using classical ones by iterating in distributions, see Lemma \ref{L2.1} below.  This   enables us to construct the strong solution. However, since the approximating SDEs depend on the initial distributions, this method does not provide other properties from existing results for classical SDEs. Fortunately, we are able to develop coupling argument  to investigate the $\W_2$-exponential convergence,     Harnack inequality and applications for  the associated nonlinear semigroup.

\

  In Section 2, we  investigate the existence, uniqueness and time-space continuity of solutions. In Section 3, we study the $\W_2$-exponential contraction of $P_t^*$,   which implies the exponential ergodicity in the time-homogenous case. In Sections 4 and 5, we use coupling by change of measures to establish    Harnack and shift Harnack inequalities and make applications.   Finally, in Section 6, we apply the main results to specific models including the  homogeneous Landau equation. These results provide pointwise estimates on the distributions, which are essentially different from existing results on $L^p$-estimates and Sobolev regularities   derived in \cite{DV1,DV2,CA} for the homogeneous Landau equation.

\section{Existence, uniqueness and time-space continuity}

It was already explained in Introduction that the distribution dependence of coefficients may cause trouble in the study of DDSDEs. To get rid of the distribution dependence,  we will iterate \eqref{E1} in distributions.
To prove the convergence of solutions to  iterating SDEs,  we make the following assumptions on the  continuity,  monotonicity and growth of coefficients.

 \beg{enumerate} \item[$(H1)$] (Continuity)  For every $t\ge 0$, $b_t$ is continuous on $\R^d\times\scr P_\theta$. Moreover, there exist   increasing $K_{\si,1}, K_{\si, 2}\in C([0,\infty);(0,\infty))$ such that
 $$\|\si_t(x,\mu)- \si_t(y,\nu)\|^2\le K_{\si,1}(t) |x-y|^2+ K_{\si,2}\W_\theta(\mu,\nu)^2,\ \ t\ge 0, x,y\in \R^d, \mu,\nu\in
\scr P_\theta.$$
 \item[$(H2)$] (Monotonicity) There exists  increasing $K_{b,1}, K_{b,2}\in C([0,\infty);(0,\infty))$ such that
\beg{align*} & 2 \big\<b_t(x,\mu)-b_t(y,\nu), x-y\big\>
  \le K_{b,1}(t) |x-y|^2+ K_{b,2}(t)\W_\theta(\mu,\nu) |x-y|,\\
  &\qquad \ \ t\ge 0; x,y\in \R^d; \mu,\nu\in
\scr P_\theta.\end{align*}
\item[$(H3)$] (Growth)  $b$ is bounded on bounded sets in $[0,\infty)\times \R^d\times \scr P_\theta$, and there exists  increasing $K_{b,3}\in C([0,\infty);(0,\infty))$ such that
$$|b_t(0,\mu)|^\theta\le K_{b,3}(t) \big\{1 +\mu(|\cdot|^\theta)\big\},\ \  x\in\R^d, t\ge 0, \mu\in \scr P_\theta.$$
\end{enumerate}

\subsection{Main results}

We first consider the existence, uniqueness and $\W_\theta$-Lipschitz continuous in initial distributions.

\beg{thm}\label{T1.1} Assume $(H1)$-$(H3)$ for some $\theta\in [1,\infty)$ such that $K_{\si,2}=0$ when $\theta<2$.
 \beg{enumerate} \item[$(1)$] The DDSDE $\eqref{E1}$ has strong/weak existence and uniqueness for every initial distribution  in $\scr P_\theta$.
   Moreover, for any $p\ge \theta$ and $s\ge 0$, $\E|X_{s,s}|^{p}<\infty$ implies
\beq\label{ES1} \E \sup_{t\in [s,T]}|X_{s,t} |^{p}<\infty,\ \ T \ge s.\end{equation}
\item[$(2)$] There exists increasing $\psi: [0,\infty)\to [0,\infty)$ such that for any two solutions $X_{s,t}$ and $Y_{s,t}$ of $\eqref{E1}$ with $\L_{X_{s,s}},\L_{Y_{s,s}}\in \scr P_\theta$,
\beq\label{CC1} \E |X_{s,t}-Y_{s,t}|^\theta\le \big(\E|X_{s,s}-Y_{s,s}|^\theta\big)\e^{ \int_s^t\psi(r)\d r},\ \ t\ge s\ge 0.\end{equation} Consequently,
\beq\label{CC2}\lim_{\E|X_{s,s}-Y_{s,s}|^\theta\to 0} \P\Big(\sup_{r\in [s,t]}|X_{s,r}-Y_{s,r}|\ge \vv\Big)=0,\ \ t>s\ge 0,\vv>0;\end{equation} and
\beq\label{X0}  \W_\theta(P_{s,t}^*\mu_0, P_{s,t}^*\nu_0)^\theta\le   \W_2(\mu_0,\nu_0)^\theta\e^{\int_s^t\psi(r)\d r},\ \ t\ge s\ge 0.\end{equation}
\end{enumerate} \end{thm}

 Since  assumptions $(H1)$-$(H3)$ are weaker for larger $\theta$,   and $\dd_x\in \scr P_\theta$ for any $\theta\ge 1$,   by Theorem \ref{T1.1} the DDSDE \eqref{E1} has a unique solution    $(X_t(x))_{t\ge 0}$ for $X_0=x$. The next result says that $X$ is continuous in $(t,x)\in [0,\infty)\times \R^d$ provided $b$ has a polynomial growth. Because  the coefficients depend  on the distribution of solution, it seems hard
to prove the flow property, for instance to prove that $\P$-a.s. for all $t$ the map $X_t(\cdot):\R^d\to \R^d$ is a diffeormorphism, by using techniques developed in the classical setting. So, we leave the study to the future.

\beg{thm}\label{T2.2} Assume $(H1)$-$(H3)$ for some $\theta\ge 1$. If there exists $p\ge 1$ such that
\beq\label{X1} |b_t(x,\mu)|\le K(t) \big\{1+|x|^p+ \mu(|\cdot|^p)\big\},\ \ t\ge 0, x\in \R^d\end{equation} holds for some increasing function
$K: [0,\infty)\to (0,\infty)$, then $\P$-a.s. the map  $$[0,\infty)\times \R^d\ni (t,x)\mapsto X_t(x)\in \R^d$$ is continuous. \end{thm}

To prove these results, we first approximate \eqref{E1} using classical SDEs by iterating in distributions.

\subsection{An approximation argument using classical SDEs}

We fixed   $s\ge 0$ and $\F_s$-measurable random variable $X_{s,s}$ on $\R^d$ with $\E|X_{s,s}|^\theta<\infty$. Let
$$X_{s,t}^{(0)}=X_{s,s},\ \ \mu_{s,t}^{(0)}=\L_{X_{s,s}},\ \ t\ge s.$$  For any $n\ge 1,$ let $(X_{s,t}^{(n)})_{t\ge s}$ solve the classical SDE
\beq\label{EN} \d X_{s,t}^{(n)}= b_t(X_{s,t}^{(n)}, \mu_{s,t}^{(n-1)}) \d t + \si_t(X_{s,t}^{(n)},\mu_{s,t}^{(n-1)})\,\d W_t,\ \ X_{s,s}^{(n)}=X_{s,s}, t\ge s,\end{equation}
where $\mu_{s,t}^{(n-1)}=\L_{X_{s,t}^{(n-1)}}.$

\beg{lem} \label{L2.1} Assume $(H1)$-$(H3)$ for some $\theta\in [1,\infty)$.
\beg{enumerate} \item[$(1)$]
For every $n\ge 1$, the SDE $\eqref{EN}$ has a unique strong solution and
\beq\label{*2} \E\sup_{t\in [s,T]} |X_{s,t}^{(n)}|^\theta<\infty,\ \ T>s, n\ge 1.\end{equation}
 \item[$(2)$] If either $\theta\ge 2$ or $\si(x,\mu)$ does not depend on $\mu$, then for any $T>0$ there exists $t_0>0$ which is independent on $s\in [0,T]$ and $X_{s,s}$, such that
\beq\label{*3} \E \sup_{ t\in [s, s+t_0]} |X_{s,t}^{(n+1)}-X_{s,t}^{(n)}|^\theta \le 2^\theta \e^{-n}  \E\sup_{t\in [s,s+t_0]} |X_{s,t}^{(1)}|^2,\ \ s\in [0,T], n\ge 1.\end{equation}\end{enumerate}
\end{lem}

\beg{proof}   Without loss of generality, we only prove for  $s=0$.

(1) We first prove that  the SDE \eqref{EN} has a unique strong solution and \eqref{*2} holds.
 By $(H1)$, $b_t(x, \mu_t^{(0)})$  and $\si_t(x, \mu_t^{(0)})$ are continuous in $x$. Then  the SDE \eqref{EN} for $n=1$  has a weak solution up to life time (see   \cite[Theorem 6.1.6]{SV} and  \cite[p.155-163]{IW}). Next, by   It\^o's formula it is easy to see that $(H2)$ implies  the pathwise uniqueness. According to the Yamada-Watanabe principle \cite{YW}, the SDE has a unique solution up to life time.
It remains to prove \eqref{*2}. By $(H3)$ and It\^o's formula we have
\beq\label{*c} \beg{split}  \d |X_t^{(1)}|^2  &=2\<\si_t(X_t^{(1)},\mu_t^{(0)})\d W_t, X_t^{(1)}\>\\
&\quad + \big\{ 2\big\<b_t(X_t^{(1)},  \mu_t^{(0)}), X_t^{(1)}\big\> + \|\si_t(X_t^{(1)},\mu_t^{(0)})\|_{HS}^2\big\}\d t.\end{split} \end{equation}
By $(H1)$ with $y=0,\nu=\dd_0$, we have
\beq\label{H1*} \|\si_t(x,\mu)\|_{HS}^2\le K(t) \big\{1+|x|^2+\mu(|\cdot|^\theta)^{\ff 2 \theta}\big\}\end{equation} for some increasing $K: [0,\infty)\to [0,\infty).$ Combining this with $(H2)$ and $(H3)$, we may find   increasing $H:[0,\infty)\to (0,\infty)$ such that
\beg{align*}& \max\big\{2\big\<b_t(x,\mu_t^{(0)}), x\big\>,\ \|\si_t(x,\mu_t^{(0)})\|_{HS}^2\big\}\\
&\le \max\{2\big\<b_t(x,\mu_t^{(0)})-b_t(0,\mu_t^{(0)}), x \big\>
  +   2|b_t(0,\mu_t^{(0)})|\cdot |x|, \   \|\si_t(x,\mu_t^{(0)})\|^2_{HS}\big\}\\
&\le H(t) \big\{1+ |x|^2+\mu_t^{(0)}(|\cdot|^\theta)^{\ff 2 \theta}\big\},\ \ \ t\ge 0, x\in \R^d.\end{align*}
Thus, by \eqref{*c}, \eqref{H1*} and using It\^o's formula, there exists a constant $c_1(\theta)>0$ such that
\beg{align*} \d (1+|X_t^{(1)}|^2)^{\ff\theta 2} \le &\theta (1+|X_t^{(1)}|^2)^{\ff{\theta-2}2} \<\si_t(X_t^{(1)},\mu_t^{(0)})\d W_t, X_t^{(1)}\>\\
&+c_1(\theta) H(t) \Big\{(1+|X_t^{(1)}|^2)^{\ff \theta 2} +\mu_t^{(0)} (|\cdot|^\theta)^{\ff 2 \theta\lor 1}\Big\}\d t.\end{align*}
Letting $\tau_N= \inf\{t\ge 0: |X_t^{(1)}|\ge N\}$, we conclude from this,   \eqref{H1*} and   the BDG inequality yield  that    for some
increasing $\Psi: [0,\infty)\to (0,\infty)$,
\beg{align*} &\E \sup_{s\in [0,t\land \tau_N]} \big(1+ |X_s^{(1)}|^2\big)^{\ff \theta 2} \le c_1(\theta)
 H(t) \E\int_0^{t\land \tau_N} \Big\{\big(1+ |X_s^{(1)}|^2\big)^{\ff \theta 2} + \mu_t^{(0)}(|\cdot|^\theta)^{\ff 2\theta\lor 1}\Big\}\d s \\
&\qquad + c_2(\theta) K(t) \E\bigg(\int_0^{t\land \tau_N} \Big\{\big(1+|X_s^{(1)}|^2\big)^\theta   + \big(1+|X_s^{(1)}|^2\big)^{\theta-\ff 1 2} \mu_s^{(0)}(|\cdot|^\theta)^{\ff 2\theta\lor 1}\Big\}\d s\bigg)^{\ff 1 2}\\
&\le \Psi(t) \E \int_0^{t\land \tau_N} \Big\{\big(1+|X_s^{(1)}|^2\big)^{\ff\theta  2}+ \mu_s^{(0)}(|\cdot|^\theta)^{\ff 2 \theta\lor 1} \Big\}\d s
  + \ff 1 2 \E \sup_{s\in [0,t\land \tau_N]}\big(1+|X_s^{(1)}|^2\big)^{\ff\theta  2}.\end{align*}
Therefore,
$$\E \sup_{s\in [0,t\land \tau_N]} \big(1+|X_s^{(1)}|^2\big)^{\ff\theta  2}\le 2\Psi(t) \int_0^{t} \Big\{\big(1+|X_{s\land\tau_N}^{(1)}|^2\big)^{\ff\theta  2} + \mu_s^{(0)}(|\cdot|^\theta)^{\ff 2 \theta\lor 1}\Big\}\d s.$$ By Gronwall's lemma and letting $N\to\infty$, we arrive at
$$  \E \sup_{s\in [0,t]} \big(1+|X_s^{(1)}|^2\big)^{\ff\theta  2}\le\Big(1+  2t\Psi(t)\sup_{s\in [0,t]}\big(\E |X_s^{(0)}|^\theta\big)^{\ff 2 \theta\lor 1} \Big) \exp\big[2t\Psi(t)\big]<\infty.$$
Therefore,   \eqref{*2}  holds  for $n=1$.

Now, assume that the assertion holds for $n=k$ for some $k\ge 1$, we intend to prove it for $n=k+1$. This can be done in the same way by using
$(X_\cdot^{(k+1)}, \mu_\cdot^{(k)}, X_\cdot^{(k)})$ in place of $(X_\cdot^{(1)}, \mu_\cdot^{(0)},X_\cdot^{(0)})$. So, we omit the proof to save space.

(2) To prove \eqref{*3}, for $n\ge 1$ we simply denote
\beg{align*} &\xi_t^{(n)} = X_t^{(n+1)}- X_t^{(n)},\\
&\LL_t^{(n)}= \si_t(X_t^{(n+1)},\mu_t^{(n)})- \si_t(X_t^{(n)},\mu_t^{(n-1)}),\\
&B_t^{(n)}= b_t(X_t^{(n+1)}, \mu_t^{(n)} ) - b_t(X_t^{(n)}, \mu_t^{(n-1)} ).\end{align*}
Below we prove for 1) $\theta\ge 2$ and 2) $\theta<2$ but $K_{\si, 2}=0$ respectively.

Let $\theta\ge 2$.  By $(H1)$, $(H2)$ and It\^o's formula,  there exists increasing $K_0: [0,\infty)\to [0,\infty)$ such that
\beq\label{*f} \d |\xi^{(n)}_t|^2  \le 2 \<\LL_t^{(n)} \d W_t, \xi_t^{(n)}\>+  K_0(t) \big\{|\xi_t^{(n)}|^2 + \W_\theta(\mu_t^{(n)}, \mu_t^{(n-1)})^2\big\}\d t.\end{equation}
Combining this with $(H1)$ and using the BDG inequality, we may find out increasing functions $K_1, K_2: [0,\infty)\to (0,\infty)$ such that
 \beg{align*}& \E \sup_{s\in [0,t]}  |\xi_s^{(n)}|^\theta\le 2^{\ff \theta 2 -1}\E \bigg|\sup_{s\in [0,t]}\int_0^s 2 \<\LL_t^{(n)} \d W_t, \xi_t^{(n)}\>\bigg|^{\ff \theta 2} \\
 &\quad + 2^{\ff \theta 2 -1} \E\bigg(\int_0^t   K_0(s) \big\{|\xi_s^{(n)}|^2 + \W_\theta(\mu_s^{(n)}, \mu_s^{(n-1)})^2\big\}\d s\bigg)^{\ff \theta 2}\\
   & \le  K_1(t) \E\bigg(\int_0^t \big\{|\xi_s^{(n)}|^2 \|\LL_s^{(n)}\|^2\big\}\d s\bigg)^{\ff \theta 4}
 + K_1(t) \int_0^t \Big\{\E   |\xi_s^{(n)}|^\theta + \W_\theta(\mu_s^{(n)}, \mu_s^{(n-1)})^\theta\Big\} \d s\\
 &\le \ff 1 2\E \sup_{s\in [0,t]}  |\xi_s^{(n)}|^\theta + K_2(t) \int_0^t \Big\{\E   |\xi_s^{(n)}|^\theta + \W_\theta(\mu_s^{(n)}, \mu_s^{(n-1)})^\theta\Big\} \d s.  \end{align*}   Then
$$\E \sup_{s\in [0,t]}  |\xi_s^{(n)}|^2 \le 2 K_2(t) \int_0^t \Big\{\E   |\xi_s^{(n)}|^\theta + \W_\theta(\mu_s^{(n)}, \mu_s^{(n-1)})^\theta\Big\} \d s,\ \ t\ge 0. $$   By
  Gronwall's lemma, we obtain
\beq\label{*g} \beg{split} \E \sup_{s\in [0,t]}  |\xi_s^{(n)}|^\theta& \le 2tK_2(t) \e^{2tK_2(t)} \sup_{s\in [0,t]} \W_\theta(\mu_s^{(n)}, \mu_s^{(n-1)})^\theta\\
&\le 2tK_2(t) \e^{2tK_2(t)}  \E  \sup_{s\in [0,t]}|\xi_s^{(n-1)}|^\theta,\ \  t\ge 0.\end{split}\end{equation}
 Taking $t_0>0$ such that $ 2 t_0K_2(t_0) \e^{2t_0K_2(t_0)}\le \e^{-1}$, we arrive  at
 $$\E \sup_{s\in [0,t_0]} |\xi_s^{(n)}|^\theta\le \e^{-1}\E \sup_{s\in [0,t_0]} |\xi_s^{(n-1)}|^\theta,\ \ n\ge 1.$$ Since
 $$\E \sup_{s\in [0,t_0]} |\xi_s^{(0)}|^\theta \le 2^{\theta-1}  \E  \Big\{|X_0|^\theta+ \sup_{s\in [0,t_0]} |X_s^{(1)}|^\theta\Big\}\le 2^\theta \E \sup_{s\in [0,t_0]}  |X_s^{(1)}|^2,$$ we prove \eqref{*3}.

 Let $\theta\in [1,2)$ but $K_{\si,2}=0$.  Then instead of \eqref{*f} we have
 $$\d |\xi^{(n)}_t|^2  \le 2 \<\LL_t^{(n)} \d W_t, \xi_t^{(n)}\>+  K_0(t) |\xi_t^{(n)}|\big\{|\xi_t^{(n)}|  + \W_\theta(\mu_t^{(n)}, \mu_t^{(n-1)})\big\}\d t.$$ Since $\theta\le 2$, for any $\vv>0$, by It\^o's formula we obtain
 $$\d (\vv+|X_t^{(n)}|^2)^{\ff \theta 2} \le \theta (\vv+|X_t^{(n)}|^2)^{\ff {\theta -2}2}\Big\{\<\LL_t^{(n)} \d W_t, \xi_t^{(n)}\>+  \ff{K_0(t) } 2 |\xi_t^{(n)}|\big\{|\xi_t^{(n)}|  + \W_\theta(\mu_t^{(n)}, \mu_t^{(n-1)})\big\}\d t\Big\}.$$
 Since $(H1)$ with $K_{\si, 2}=0$ implies $\|\LL_t^{(n)}\|^2\le K_{\si, 1}(t)|\xi_t^{(n)}|^2,$  this and  the BDG inequality yield
 \beg{align*} \E \sup_{s\in [0,t]} (\vv+|X_s^{(n)}|^2)^{\ff \theta 2}
 &\le  K_1(t)  \E \bigg(\int_0^t  (\vv+|X_s^{(n)}|^2)^{\theta}\d s \bigg)^{\ff 1 2}\\
  &\quad + K_1(t) \E \int_0^t \Big\{ (\vv+|X_s^{(n)}|^2)^{\ff \theta 2} + (\vv+|X_s^{(n)}|^2)^{\ff {\theta-1} 2}W_\theta (\mu_s^{(n)}, \mu_s^{(n-1)})\Big\}\d s\\
  &\le \ff 1 2 \E \sup_{s\in [0,t]} (\vv+|X_s^{(n)}|^2)^{\ff \theta 2}+ K_2(t)\int_0^t \Big\{ (\vv+|X_s^{(n)}|^2)^{\ff \theta 2}
  +W_\theta(\mu_s^{(n)}, \mu_s^{(n-1)})^\theta\Big\}\d s \end{align*} for some increasing $K_1, K_2: [0,\infty)\to [0,\infty)$. Letting $\vv\to 0$ and using Gronwall's inequality, we prove \eqref{*g}, which implies the desired estimate \eqref{*3}  as explained above.
\end{proof}

\subsection{Proofs of Theorem \ref{T1.1} and Theorem \ref{T2.2}}

\beg{proof}[Proof of Theorem \ref{T1.1}] Without loss of generality, we only consider the DDSDE \eqref{E1} from time $s=0$.

(1) Since the uniqueness follows from \eqref{CC2}  which will be proved in the next step, in this step we only  prove the existence and the estimate  \eqref{ES1}.

By Lemma \ref{L2.1}, there exists an adapted continuous process $(X_t)_{t\in [0,t_0]}$ such that
\beq\label{A01}\lim_{n\to\infty} \sup_{t\in [0,t_0]} \W_\theta(\mu_t^{(n)},\mu_t)^\theta\le \lim_{n\to\infty} \E \sup_{t\in [0,t_0]}  |X_t^{(n)}- X_t|^\theta=0,\end{equation} where
 $\mu_t$ is the distribution of $X_t$. Noting that due to \eqref{EN}
$$X_t^{(n)}=X_0+ \int_0^t  b_s(X_s^{(n)},\mu_s^{(n-1)}) \d s +\int_0^t\si_s(X_s^{(n)},\mu_s^{(n-1)})\d W_s, $$ it follows from \eqref{A01}, $(H1)$ and  $(H3)$  that $\P$-a.s.
$$X_t= X_0+\int_0^t b_s(X_s, \mu_s)\d s +\int_0^t \si_s(X_s, \mu_s)\d W_s,\ \ t\in [0,t_0].$$
Therefore, $(X_t)_{t\in [0,t_0]}$ is a solution to \eqref{E1}, and \eqref{A01} implies  $ \E \sup_{s\in [0,t_0]} |X_s|^\theta<\infty.$
   Since $t_0>0$ is independent of $X_0$,  we conclude that  \eqref{E1} has a unique solution
  $(X_t)_{t\ge 0}$ with
  \beq\label{*X} \E \sup_{s\in [0,t]}|X_s|^\theta<\infty,\ \ \ t\in (0,\infty).\end{equation}

It remains to prove \eqref{ES1} for $\E|X_0|^{p}<\infty$. As in the proof of \eqref{*2} above, by $(H1)$-$(H3)$ and It\^o's formula we have
  $$\d |X_t|^2 \le 2 \<\si_t(X_t, \L_{X_t})\d W_t, X_t\> + H(t) \big\{1+|X_t|^2+(\E |X_t|^\theta)^{\ff 2 \theta}\big\}\d t $$   for some increasing function $H: [0,\infty)\to (0,\infty).$ Then applying It\^o's formula to $(1+|X_t|^2)^{\ff p 2}$ and repeating step (1) in the proof of Lemma \ref{L2.1}, we prove  \eqref{ES1}.

(2) By It\^o's formula, $(H2)$ and $(H1)$ with $K_{\si,2}=0$ if $\theta<2$, we have
\beq\label{cc} \beg{split} &\d |X_t-Y_t|^2\le   2\big\<X_t-Y_t, \{\si_t(X_t,\L_{X_t})-\si_t(Y_t,\L_{Y_t})\}\d W_t\big\>\\
&+K_1(t)\big\{|X_t-Y_t|^2+1_{\{\theta\ge 2\}}\W_\theta(\L_{X_t},\L_{Y_t})^2+ |X_t-Y_t|W_\theta(\L_{X_t},\L_{Y_t}) \big\}\d t.\end{split} \end{equation}  By It\^o's formula we obtain
\beg{align*} \d |X_t-Y_t|^\theta &\le \theta |X_t-Y_t|^{\theta -2} \big\<X_t-Y_t, \{\si_t(X_t,\L_{X_t})-\si_t(Y_t,\L_{Y_t})\}\d W_t\big\>\\
&\quad + K_2(t) \big\{|X_t-Y_t|^\theta +   \W_\theta(\L_{X_t},\L_{Y_t})^\theta\big\}\d t\end{align*} for some increasing $K_2: [0,\infty)\to [0,\infty).$ Noting that $\W_\theta(\L_{X_t},\L_{Y_t})^\theta \le \E|X_t-Y_t|^\theta<\infty$, this implies
$$ \E|X_{t}-Y_{t}|^\theta \le \E|X_0-Y_0|^\theta+2 \int_0^{t} K_2(s) \E|X_{s}-Y_{s}|^2 \d s. $$  By Gronwall's lemma, we prove  \eqref{CC1}.

To prove \eqref{CC2},  let $\tau_\vv:= \inf\{t\ge 0: |X_t-Y_t|\ge \vv\}$ for $\vv\in (0,1)$. By \eqref{CC1} and \eqref{cc}, there exists increasing $K: [0,\infty)\to [0,\infty)$ such that
\beg{align*} &\E|X_{t\land \tau_\vv}-Y_{t\land\tau_\vv}|^\theta\\
&\le \E|X_0-Y_0|^\theta +\E \int_0^{t\land \tau_\vv} K(s)\big(\E|X_s-Y_s|^\theta+  |X_{s}- Y_{s}|^\theta\big)\d s\\
&\le \big\{1+t\e^{t\psi(t)}\big\}\E|X_0-Y_0|^2+ \int_0^t  K(s)\big\{ \E|X_{s\land \tau_\vv}-Y_{s\land\tau_\vv}|^\theta\big\}\d s,\ \ \ t\ge 0. \end{align*}   By Gronwall's lemma, there exists positive $\phi\in C([0,\infty))$ such that
$$\vv^\theta\P(\tau_\vv\le t)\le \E|X_{t\land \tau_\vv}-Y_{t\land\tau_\vv}|^\theta\le \phi(t) \E|X_0-Y_0|^\theta,\ \ t\ge 0.$$
Therefore,
$$\P\Big(\sup_{s\in [0,t]}|X_s-Y_s|\ge \vv\Big)= \P(\tau_\vv\le t) \le \vv^{-\theta}\psi(t) \E|X_0-Y_0|^\theta,\ \ t,\vv>0.$$
Hence,  \eqref{CC2} holds.

(3)    Let   $(X_t,W_t)$ and $(\tt X_t,\tt W_t)$  with respect to  $(\OO, \{\F_t\}_{t\ge 0}, \P)$ and
$(\tt\OO, \{\tt\F_t\}_{t\ge 0}, \tt\P)$ respectively be two weak solutions such that  $\L_{X_0}|_{\P}= \L_{\tt X_0}|_{\tt\P}$.
Then $X_t$ solves \eqref{E1} while  $\tt X_t$ solves
\beq\label{E1'} \d \tt X_t = b_t(\tt X_t, \L_{\tt X_t}|_{\tt\P})\d t + \si_t(\tt X_t, \L_{\tt X_t}|_{\tt\P})\d \tt W_t.\end{equation}
To prove that $\L_{X}|_{\P}=\L_{\tt X}|_{\tt\P}$, let $\mu_t= \L_{X_t}|_{\P}$ and
$$\bar b_t(x)= b_t(x, \mu_t),\ \ \bar \si_t(x)= \si(x,\mu_t),\ \ t\ge 0, x\in\R^d.$$
  By $(H1)$-$(H3)$, the   SDE
\beq\label{E10} \d \bar X_t = \bar b_t(\bar X_t)\d t + \bar \si_t(\bar X_t)\d \tt W_t,\ \ \bar X_0= \tt X_0 \end{equation}
has a unique solution for any initial points.
According to Yamada-Watanabe, it also has weak uniqueness. Noting that
$$\d X_t = \bar b_t( X_t)\d t + \bar \si_t(  X_t)\d  W_t,\ \ \L_{X_0}|_{\P}= \L_{\tt X_0}|_{\tt\P},$$ the weak uniqueness of \eqref{E10} implies
\beq\label{HW} \L_{\bar X}|_{\tt\P}= \L_X|_{\P}.\end{equation}
So, \eqref{E10} reduces to
$$ \d \bar X_t = b_t(\bar X_t, \L_{\bar X_t}|_{\tt\P})\d t + \si_t(\bar X_t, \L_{\bar X_t}|_{\tt\P})\d \tt W_t,\ \ \bar X_0=\tt X_0.$$
Since by (1)  the  DDSDE \eqref{E1'} has a unique solution, we obtain  $\bar X=\tt X$. Therefore, the  weak uniqueness follows from \eqref{HW}.

Finally, for any $\mu_0,\nu_0\in \scr P_\theta$,  take $\F_0$-measurable random variables $X_0,Y_0$ such that $\L_{X_0}=\mu_0, \L_{Y_0}=\nu_0$ and
$\W_\theta(\mu_0, \nu_0)^\theta=\E|X_0-Y_0|^\theta$. Since $\W_\theta(P_t^*\mu_0, P_t^*\nu_0)^\theta\le \E|X_t-Y_t|^\theta,$  \eqref{CC1} implies \eqref{X0}.

\end{proof}

\beg{proof}[Proof of Theorem \ref{T2.2}] Since the assumptions are weaker for larger $\theta$, we may and do assume that $\theta\ge 2.$ By Kolmogorov's modification theorem, it suffices to prove
\beq\label{X4} \E |X_t(x)-X_s(y)|^m\le \Phi(s,t;x,y)(|x-y| + |s-t|)^q,\ \ |t-s|+|x-y|\le 1\end{equation}
for some constants $m>0, q>1$ and locally bounded function $\Phi$ on $[0,\infty)^2\times\R^{2d}.$   Firstly, by \eqref{ES1} and \eqref{CC1}, we may find out an increasing function $\psi: [0,\infty)\to (0,\infty)$ such that
\beq\label{X5} \beg{split} \E |X_t(x)-X_t(y)|^{2\theta} &\le \big\{\E|X_t(x)-X_t(y)|^\theta\big\}^{\ff 2 3} \times\big\{\E|X_t(x)-X_t(y)|^{4\theta}\big\}^{\ff 1 3}\\
&\le \psi(t)(1+|x|+|y|)^{\ff {4\theta} 3}|x-y|^{\ff {2\theta} 3},\ \ t\ge 0, x,y\in \R^d.\end{split}\end{equation}
Next, by $(H3)$ and \eqref{X1}, there exist a constant $C>0$ and an increasing function $\phi: [0,\infty)\to (0,\infty)$ such that
\beg{align*} \E|X_t(x)-X_s(x)|^{2\theta} &\le C \E\bigg(\int_s^t K(r)(1+|X_r(x)|^p+\E |X_r(x)|^p)\d r\bigg)^{2\theta}\\
&\quad +C \E\bigg(\int_s^t K_2(r) (1+|X_r(x)|^2+ \E |X_r(x)|^2)\d r\bigg)^\theta\\
&\le \phi(t) (1+|x|^{2\theta p})(t-s)^\theta,\ \ |t-s|\le 1, x\in \R^d.\end{align*} This together with \eqref{X5} implies the desired \eqref{X4} with $p=4, q= \ff {2\theta} 3>1.$
\end{proof}

\section{$\W_2$-Exponential  contraction of $P_{s,t}^*$}

We intend to estimate  the Wasserstein distance of solutions with different initial distributions and investigate the exponential ergodicity. For simplicity, we only consider the  $\W_2$-distance.   To this end, we use the following condition to replace $(H2)$:
\beg{enumerate}\item[$(H2')$] There exist   positive functions $C_2,C_2\in L^1_{loc}(\d t)$ such that
\beg{align*} & 2\<b_t(x,\mu)- b_t(y,\nu), x-y\> + \|\si_t(x,\mu)-\si_t(y,\nu)\|_{HS}^2\\
&\le C_1(t) \W_2(\mu,\nu)^2 - C_2(t) |x-y|^2,\ \
 t\ge 0; x,y\in \R^d; \mu,\nu\in\scr P_2.\end{align*}
 \end{enumerate}

In the following result we present estimates on Wasserstein distance of $P_{s,t}^*$, which in particular provide  exponential upper bound estimates for homogeneous Landuan equation with Maxwell molecules, see Corollary  \ref{CC} below. See also \cite{Ta} for Wasserstein contraction of Boltzmann equation with Maxwell molecules without estimates on the convergence rate, where   SDEs are driven by Poisson point processes are applied. 

\beg{thm}\label{T1.2} Assume $(H1)$, $(H2')$ and $(H3)$.
\beg{enumerate} \item[$(1)$] For any $\mu_0,\nu_0\in \scr P_2,$
$$\W_2(P_{s,t}^*\mu_0,P_{s,t}^*\nu_0)^2\le \W_2(\mu_0,\nu_0)^2\e^{\int_s^t\{C_1(r)-C_2(r)\}\d r},\ \ t\ge 0.$$
\item[$(2)$] Let $b_t=b$ and $\si_t=\si$ do not depend on time $t$ such that $(H2')$ holds for some constants $C_1$ and $C_2$.    If $C_2>C_1$ then $P_t$ has a   unique invariant probability measure $\mu\in\scr P_2$ such that
$$\W_2(P_t^*\nu_0,\mu)^2\le \W_2(\nu_0,\mu)^2\e^{-(C_2-C_1)t},\ \ t\ge 0, \nu_0\in \scr P_2.$$
\end{enumerate}
\end{thm}

\beg{proof}  (1) Without loss of generality, we only prove for $s=0$. Let $X_t$ and $Y_t$ be two solutions to \eqref{E1} such that $\L_{X_0}=\mu_0, \L_{Y_0}=\nu_0$ and
\beq\label{X2} \W_2(\mu_0,\nu_0)^2=\E|X_0-Y_0|^2.\end{equation}
Simply denote $\mu_t=\L_{X_t}, \nu_t=\L_{Y_t}, t\ge 0.$   By $(H2')$ and It\^o's formula we have
\beg{align*} \d |X_t-Y_t |^2 &\le 2 \big\<X_t -Y_t, \{\si_t(X_t,\mu_t)- \si_t(Y_t,\nu_t)\} \d W_t\big\>\\
&\quad +\big\{C_1(t)  \W_2(\mu_t,\nu_t)^2-C_2(t) |X_t-Y_t|^2  \big\}\d t.\end{align*} Noting that $\W_2(\mu_s,\nu_s)^2\le \E|X_s-Y_s|^2,$
combining this with \eqref{X2}
we obtain
$$\E |X_t -Y_t|^2\le \W_2(\mu_0,\nu_0)^2 +   \int_0^t\big\{[C_1(s)-C_2(s)] \E|X_s -Y_s|^2\big\}\d s.$$ This implies  the first assertion    by  Gronwall's lemma.

(2)  Let $\dd_0$ be the Dirac measure at point $0\in\R^d$. Then $P_t^*\dd_0= \L_{X_t(0)}$.   We first prove
\beq\label{CV} \lim_{t\to\infty} \W_2(P_t^*\dd_0,\mu)=0\end{equation} for some $\mu\in \scr P_2.$
To this end, it suffices to show that $\{P_t^*\dd_0\}_{t\ge 0}$ is a $\W_2$-Cauchy family when $t\to\infty$; that is,
\beq\label{CV2}\lim_{t\to\infty} \sup_{s\ge 0} \W_2(P_t^*\dd_0, P_{t+s}^*\dd_0)=0.\end{equation}
We will prove this using the shift-coupling and the weak uniqueness according to Theorem \ref{T1.1}(3). More precisely, for any $s\ge 0$, it is easy to see that $(\bar X_t:=X_{t+s}(0))_{t\ge 0}$ solves the DDSDE
$$ \d \bar X_{t}= b(\bar X_{t}, \L_{\bar X_{t}})\d t+ \si(\bar X_{t}, \L_{\bar X_{t}})\d \bar W_{t},\ \ \bar X_{0}= X_s(0)$$ for the $d$-dimensional Brownian motion $\bar W_t:= W_{t+s}- W_s.$ So, by the weak uniqueness we have
\beq\label{AL} P_t^*(P_s^*\dd_0)= \scr L_{\bar X_t}= \scr L_{X_{t+s}(0)} = P_{t+s}^*\dd_0,\ \ s,t\ge 0.\end{equation}   Combining this with Theorem \ref{T1.2}(1) and letting $X_t(P_s^*\dd_0)$ solve \eqref{E1} with $\L_{X_0}=P_s^*\dd_0$,  we obtain
\beg{align*}& \W_2(P_{t+s}^*\dd_0, P_t^*\dd_0)^2= \W_2(\L_{X_t(P_s^*\dd_0)}, \L_{X_t(0)} )^2\\
&\le \W_2(P_s^*\dd_0,\dd_0)^2 \e^{-(C_2-C_1)t}= \e^{-(C_2-C_1)t} \E|X_s(0)|^2,\ \ s,t\ge 0.\end{align*}
This implies \eqref{CV2} provided
\beq\label{CV3}\sup_{s\ge 0} \E |X_s(0)|^2<\infty. \end{equation}
By $(H2')$ and $(H3)$ for constant $C_1<C_2$ and $K_2$, it is easy to see that
$$2\<b(x,\mu),x\>+\|\si(x,\mu)\|_{HS}^2\le C_0 -(C_2-\vv)|x|^2+ (C_1+\vv)\mu(|\cdot|^2)$$
holds for some constant $C_0>0$ and $\vv:= \ff{C_2-C_1}4>0.$
By It\^o's formula and Gronwall's lemma, this implies
$$\E|X_t(0)|^2\le C_0\e^{-(C_2-C_1-2\vv)t},\ \ t\ge 0.$$ Therefore, \eqref{CV3} holds.

Moreover, by \eqref{X0} and \eqref{CV} we have
$$\lim_{t\to \infty} \W_2(P_s^* \mu, P_s^*P_t^*\dd_0)=0,\ \ s\ge 0.$$ Combining this with
 \eqref{CV} and \eqref{AL},     we obtain
$$\W_2(P_s^*\mu,\mu)\le  \lim_{t\to\infty} \W_2(P_s^*P_t^*\dd_0, P_t^*\dd_0) = \lim_{t\to\infty} \W_2(P_{t+s}^*\dd_0,P_t^*\dd_0)  = 0.$$ Then $\mu$ is an invariant probability measure.
Therefore, by  Theorem \ref{T1.2}(1) with $C_2>C_1$, for any $\nu_0\in \scr P_2$ we have
$$ \W_2(P_t^*\nu_0, \mu)^2=   \W_2(P_t^*\nu_0,P_t^*\mu)^2 \le \e^{-(C_2-C_1)t} \W_2(\nu_0,\mu)^2,\ \ t\ge 0,$$ so that the proof is finished.
 \end{proof}

\section{Harnack inequality and applications}

In this section, we  investigate the dimension-free  Harnack inequality  in the sense of \cite{W97} and the log-Harnack inequality introduced in \cite{RW10,W10} for the  DDSDE \eqref{E1}, see \cite{Wbook} and references within for general results on these type Harnack inequalities and applications.  To establish Harnack inequalities for DDSDEs  using  coupling by change of measures, we need to assume that the noise part is distribution-free; that is, we consider the following special version of \eqref{E1}:
\beq\label{E11} \d X_t= b_t(X_t,\L_{X_t})\d t +\si_t(X_t)\d W_t.\end{equation}  Then
$$(P_tf)(\mu_0)= \int_{\R^d} f  \d(P_t^*\mu_0)= \E f(X_t(\mu_0)),\ \ f\in \B_b(\R^d), t\ge 0, \mu_0\in \scr P_2,$$ where $X_t(\mu_0)$ solves \eqref{E11} with
initial distribution $\mu_0.$

To make the study  easy to follow, we first introduce the main steps in establishing Harnack inequalities using coupling by change of measures summarized in
\cite[\S 1.1]{Wbook}.
\beg{enumerate}
\item[(S1)] Let $(X_t)_{t\ge 0}$ solve \eqref{E11} with $\L_{X_0}=\mu_0$. By the uniqueness we have   $\mu_t:=\P_t^*\mu_0= \L_{X_t}$, and  the equation \eqref{E11} reduces to
\beq\label{CP1} \d X_t= b_t(X_t,\mu_t)\d t+\si_t(X_t)\d W_t.\end{equation}
\item[(S2)] Construct a process $(Y_t)_{t\in [0,T]}$ such that for a weighted probability measure $\Q:=R_T\P,$
\beq\label{CP2} X_T=Y_T\ \Q\text{-a.s., \ \ and}\ \L_{Y_T}|_\Q=P_T^*\nu_0=:\nu_T. \end{equation}
\end{enumerate}

Obviously,  (S1) and (S2) imply
\beq\label{CP0} (P_Tf)(\mu_0)= \E[f(X_T)]\ \text{and}\  (P_T f)(\nu_0)= \E_\Q[f(Y_T)]=\E[R_Tf(X_T)],\ \ f\in \B_b(\R^d).\end{equation}
Combining this with  \eqref{CP0} and Young's inequality (see \cite[Lemma 2.4]{ATW}), we obtain the log-Harnack inequality:
\beq\label{LHI}\beg{split} (P_T\log f)(\nu_0)&\le \E[R_T\log R_T]+ \log\E[f(X_T)]\\
&=\log (P_Tf)(\mu_0)+ \E[R_T\log R_T],\ \ f\in \B_b^+(\R^d);\end{split}\end{equation}
while using H\"older's inequality we prove the Harnack inequality with power $p>1$:
\beq\label{HI} (P_Tf(\nu_0))^p= (\E[R_T f(X_T)])^p\le (\E R_T^{\ff p{p-1}})^{p-1} (P_Tf^p)(\mu_0),\ \ f\in \B_b^+(\R^d).\end{equation}

To construct  $Y_t$ in (S2),    we will need the following assumption.

\beg{enumerate} \item[{\bf(A)}] $\si_t(x)$ is invertible and locally Lipschitzian in $x$ which is locally uniformly in $t\ge 0$,   and there exist increasing functions $\kk_0,\kk_1,\kk_2,\ll: [0,\infty)\to (0,\infty)$ such that for any $t\in [0,T], x,y\in\R^d$ and $\mu,\nu\in \scr P_2$, we have
\beq\label{A1}\|\si_t^{-1}\|_\infty\le \ll(t),\ \ |b_t(0,\mu)|^2+\|\si_t(x)\|^2\le \kk_0(t)(1+|x|^2+\mu(|\cdot|^2)),\end{equation}
  \beq\beg{split}\label{A2}& 2\<b_t(x,\mu)-b_t(y,\nu), x-y\> +\|\si_t(x)-\si_t(y)\|_{HS}^2 \\
  &\le \kk_1(t)|x-y|^2 +\kk_2(t) |x-y|\W_2(\mu,\nu).\end{split}
 \end{equation}\end{enumerate}

Obviously, {\bf (A)} implies assumptions $(H1)$-$(H3)$ in Theorem \ref{T1.1}.

\subsection{Main results}

 For any $\mu_0\in \scr P_2$ and $r\ge 0$, let $B(\mu_0, r)=\{\nu\in \scr P_2: \W_2(\mu_0,\nu)\le r\}.$ Let
 $$\phi (s,t) = \ll(t)^2 \bigg(\ff{\kk_1(t)}{ 1-\e^{-\kk_1(t)(t-s)} }+ \ff{t\kk_2(t)^2\exp[2(t-s)(\kk_1(t)+\kk_2(t))]}{2 }\bigg),\ \ 0\le s<t.$$
 Under  assumption {\bf (A)}, we have the following result
 for the log-Harnack inequality and regularity estimates on $P_t$.

\beg{thm}\label{T3.1} Assume {\bf (A)} and let $t>s\ge 0$.
 \beg{enumerate} \item[$(1)$] For any   $\mu_0,\nu_0\in \scr P_2$,
\beq\label{LH2}(P_{s,t}\log f)(\nu_0)\le  \log (P_{s,t}f)(\mu_0)+ \phi (s,t)\W_2(\mu_0,\nu_0)^2,\ \ f\in \B_b^+(\R^d).\end{equation}
Consequently,
\beq\label{LHH0} |\nn P_{s,t}f|^2\le 2\phi(s,t) \big\{P_{s,t} f^2 -(P_{s,t}f)^2\big\},\  \ f\in \B_b(\R^d).\end{equation}
\item[$(2)$]  For any different $\mu_0,\nu_0\in \scr P_2$, and any $f\in \B_b(\R^d)$,
\beq\label{LHH1}\beg{split}& \ff{|(P_{s,t}f)(\mu_0)-(P_{s,t}f)(\nu_0)|^2}{\W_2(\mu_0,\nu_0)^2} \\
&\le 2 \phi (s,t)\sup_{\nu\in B(\mu_0, \W_2(\mu_0,\nu_0))}\big\{(P_{s,t} f^2)(\nu)- (P_{s,t} f)^2(\nu)\big\}.\end{split}\end{equation}
Consequently,
\beq\label{LHH2}\beg{split}\|P_{s,t}^*\mu_0-P_{s,t}^*\nu_0\|_{var}&:=2\sup_{A\in\B_b(\R^d)} |(P_{s,t}^*\mu_0)(A)- (P_{s,t}^*\nu_0)(A)|\\
 &\le \ss{2\phi (s,t)} \W_2(\mu_0,\nu_0).\end{split}\end{equation}\end{enumerate}\end{thm}

Next, when $\|\si_t\|_\infty$ is locally bounded in $t$, we have the following result on Harnack inequality with power $p>1$ and applications.

\beg{thm}\label{T3.2} Assume {\bf (A)} and that   for some increasing $\gg: [0,\infty)\to (0,\infty)$,
\beq\label{VD} |\{\si_t(x)-\si_t(y)\}^*(x-y)|\le  \gg(t)|x-y|,\ \ t\ge 0.\end{equation}  Let $$p(t)= (1+4\ll(t)\gg(t))^2,\ \
 \GG(t)=  \kk_2(t)^2\ll(t)^2T\e^{2\kk_1(t)+2\kk_2(t)}.$$ Then for any    $ \mu_0,\nu_0\in \scr P_2 $ and $\F_0$-measurable random variables
$X_0, Y_0$ with $\L_{X_0}=\mu_0, \L_{Y_0}=\nu_0$,
\beq\label{HI11} \beg{split} (P_{s,t}f)^p(\nu_0)\le &(P_{s,t}f^p)(\mu_0) \exp\bigg[\ff{\ss p\,\GG(t) \W_2(\mu_0,\nu_0)^2}{(\ss p+1)[2(\ss p-1)^2-16\ll(t)^2\gg(t)^2]}\bigg]\\
&\times \bigg(\E\exp\bigg[\ff{2\ll(t)^2\kk_1(t)|X_0-Y_0|^2}{(\ss p-1)^2 (1-\e^{-\kk_1(t)(t-s)})}\bigg]\bigg)^{\ff{ \ss p\,(\ss p-1)^2 }{(\ss p+1)[2(\ss p-1)^2-16\ll(t)^2\gg(t)^2]}},\\
&\ \ \ \ \ t>s\ge 0, p\ge p(t), f\in \B_b^+(\R^d).\end{split}\end{equation}
In particular, for any $x,y\in \R^d$, $t>s\ge 0$ and $p\ge p(t)$,
\beq\label{HI2}  (P_{s,t}f)^p(x)\le  (P_{s,t}f^p)(y) \exp\bigg[\ff{\ss p  |x-y|^2\big(\GG(t)
+\ff{2\kk_1(t) \ll(t)^2 }{1-\exp[-\kk_1(t)(t-s)]}\big)}{(\ss p+1)[2(\ss p-1)^2-16\ll(t)^2\gg(t)^2]}\bigg]. \end{equation}
\end{thm}

Below we present some  consequences of the above Harnack inequalities, which provide upper bound estimates for the relative entropy between weak solutions with different initial distributions. These are essentially different from the entropy contraction of solutions derived in \cite[Theorem 1]{DV1} for solutions to \eqref{LD} with $a$ in \eqref{LLD} and $\gg \in [0,1]$:
$$\int_{\R^d} (f_t\log f_t)(x)\d x\le  \int_{\R^d} (f_0\log f_0)(x)\d x,\ \ t\ge 0.$$ Since assumption {\bf (A)} requires $\si_0$ to be invertible, the following result does not apply to the Landau equation \eqref{LD} for $a$ in \eqref{LLD}.

\beg{cor}\label{C3.3} Assume {\bf (A)} and let $t>s\ge 0$. \beg{enumerate} \item[$(1)$] For any $\mu_0,\nu_0\in\scr P_2$, $P_{s,t}^*\mu_0$ and $P_{s,t}^*\nu_0$ are equivalent and the Radon-Nykodim derivative satisfies the entropy estimate
\beq\label{BB1} \int_{\R^d} \bigg\{\log \ff{\d P_{s,t}^*\nu_0}{\d P_{s,t}^*\mu_0}\bigg\}\d P_{s,t}^*\nu_0\le \phi (s,t) \W_2(\mu_0,\nu_0)^2. \end{equation}
Consequently, in the situation of Theorem $\ref{T1.2}(2)$,
$$\int_{\R^d} \bigg\{\log \ff{\d P_{t}^*\nu_0}{\d P_{t}^*\mu}\bigg\}\d P_{t}^*\nu_0\le \phi (0, 1) \e^{-(C_2-C_1)(t-1)}\W_2(\mu,\nu_0)^2,\ \ t\ge 1.$$
\item[$(2)$] If $\eqref{VD}$ holds, then for any  $t>s\ge 0$ and $p\ge p(t)$,
\beq\label{BB2}\beg{split} & \int_{\R^d}  \bigg\{  \ff{\d P_{s,t}^*\nu_0}{\d P_{s,t}^*\mu_0}\bigg\}^{\ff 1 p} \d (P_{s,t}^*\nu_0)\\
 \le & \exp\bigg[\ff{\GG(t) \W_2(\mu_0,\nu_0)^2}{( 1+ p^{-\ff 1 2})[2(\ss p-1)^2-16\ll(t)^2\gg(t)^2]}\bigg]\\
&\times \bigg(\E\exp\bigg[\ff{2\ll(t)^2\kk_1(t)|X_0-Y_0|^2}{(\ss p-1)^2 (1-\e^{-\kk_1(t)(t-s)})}\bigg]\bigg)^{\ff{  (\ss p-1)^2 }{(1+p^{-\ff 1 2})[2(\ss p-1)^2-16\ll(t)^2\gg(t)^2]}}\end{split} \end{equation}  for   $\F_0$-measurable random variables
$X_0, Y_0$ with $\L_{X_0}=\mu_0, \L_{Y_0}=\nu_0$.
\end{enumerate}
\end{cor}

\beg{proof} According to the proof of \cite[Theorem 1.4.1]{Wbook}, when $\mu_0$ and $\nu_0$ are Dirac measures, these results follow from \eqref{LH2} and \eqref{HI2} respectively. In general,  the proof is completely similar. Precisely, for a $(P_{s,t}^*\mu_0)$-null set $A$ and $n\ge 1$, we apply \eqref{LH2} to
$f:= n 1_A+1$, so that
$$(P_{s,t}^*\nu_0)(A)\log(n+1)= (P_{s,t}^*\log f)(\nu_0)\le \phi (s,t)\W_2(\mu_0,\nu_0),\ \ n\ge 1.$$ Letting $n\to\infty$ we obtain $(P_{s,t}^*\nu_0)(A)=0,$ so that
$P_{s,t}^*\nu_0$ is absolutely continuous with respect to $P_{s,t}^*\mu_0$. By the symmetry, $P_{s,t}^*\mu_0$ is also absolutely continuous with respect to $P_{s,t}^*\nu_0$.
Moreover,   \eqref{BB1} follows from \eqref{LH2} by taking
$f=\ff{\d P_{s,t}^*\nu_0}{\d P_{s,t}^*\mu_0}$, while \eqref{BB2}  follows from \eqref{HI11} by taking
$f= (\ff{\d P_{s,t}^*\nu_0}{\d P_{s,t}^*\mu_0})^{\ff 1 p}.$ \end{proof}

\subsection{Proof of Theorem \ref{T3.1} }

Without loss of generality, we only prove for $s=0$.
As in \cite[\S2]{W11},   for fixed $T>0$,   let
\beq\label{Xi} \xi_t= \ff{1}{\kk_1(T)}\Big(1-\e^{\kk_1(T)(t-T)}\Big),\ \ t\in [0,T].\end{equation}
Let $\nu_t=P_t^*\nu_0$ and let $Y_0$ be $\F_0$-measurable with $\L_{Y_0}=\nu_0$. Consider the  SDE
\beq\label{CY} \d Y_t = \Big\{b_t(Y_t,\nu_t)+ \ff 1 {\xi_t} \si_t(Y_t)\si_t(X_t)^{-1}(X_t-Y_t)\Big\}\d t+ \si_t(Y_t) \d W_t.\end{equation}
By {\bf (A)} and $\sup_{t\in [0,T]}\nu_t(|\cdot|^2)<\infty$ due to Theorem \ref{T1.1}, this SDE has a unique solution $(Y_t)_{t\in [0,T)}$. Let
$$\tau_n:=T\land  \inf\{t\in [0,T): |X_t|+|Y_t|\ge n\},\ \ n\ge 1.$$ We have $\tau_n\uparrow T$ as $n\uparrow\infty$.
To verify step (S2), we first prove that
\beq\label{R} R_s:= \exp\bigg[\int_0^s \ff 1 {\xi_t} \big\<\si_t(X_t)^{-1}(Y_t-X_t),\d W_t\big\>-\ff 1 2 \int_0^s \ff {|\si_t(X_t)^{-1}(Y_t-X_t)|^2} {\xi_t^2}\d t\bigg]\end{equation}is a uniformly integrable martingale  for $s\in [0,T]$.

\beg{lem}\label{L3.1} Assume {\bf (A)}. Let $X_0,Y_0$ be two $\F_0$-measurable random variables such that $\L_{X_0}=\mu_0, \L_{Y_0}=\nu_0$  and    \beq\label{I1} \E|X_0-Y_0|^2= \W_2(\mu_0,\nu_0)^2.\end{equation}  Then $(R_s)_{s\in [0,T]}$ is a uniformly integrable martingale  with
\beq\label{ES} \sup_{t\in [0,T]}\E[R_t\log R_t]\le \phi (0,T)\W_2(\mu_0,\nu_0)^2.\end{equation}
\end{lem}

\beg{proof} By {\bf (A)}, for any $n\ge 1$ the process $(R_{s\land \tau_n})_{s\in [0,T)}$ is a uniformly integrable continuous martingale. Since
$\tau_n\uparrow T$ as $n\uparrow\infty$, by the martingale convergence theorem, it suffices to prove
\beq\label{ES'} \sup_{t\in [0,T],n\ge 1}\E[R_{t\land\tau_n}\log R_{t\land\tau_n}]\le \phi (0,T)  \W_2(\mu_0,\nu_0)^2.\end{equation}
We fix  $t\in (0,T)$ and $n\ge 1$.  By Girsnaov's theorem,
$$\tt W_s:= W_s -\ff 1 {\xi_s}  \si_s(X_s)^{-1}(Y_s-X_s),\ \ s\in [0, t\land \tau_n]$$   is a $d$-dimensional Brownian motion under    the weighted  probability $\Q_{t,n}:= R_{t\land \tau_n}\P$.   Reformulating \eqref{CP1} and \eqref{CY} as
\beg{equation*}\beg{split} &\d X_s=\Big\{ b_s(X_s,\mu_s) - \ff{X_s-Y_s}{\xi_s}\Big\}\d s +\si_s(X_s)\d \tt W_s,\\
&\d Y_s=  b_s(Y_s,\nu_s)\d s+   \si_s(Y_s) \d \tt W_s,\ \ s\in [0, t\land \tau_n],\end{split}\end{equation*}
by {\bf (A)} and It\^o's formula under probability $\Q_{t,n}$, we obtain
$$ \d |X_s-Y_s|^2\le \Big\{\kk_1(s)|X_s-Y_s|^2 +\kk_2(s) |X_s-Y_s|\W_2(\mu_s,\nu_s)-\ff {2|X_s-Y_s|^2}{\xi_s}\Big\}\d s +\d M_s$$ for $s\in [0, t\land \tau_n]$ and some $\Q_{t,n}$-martingale $M_s$. Then
\beq\label{V1}\beg{split}\d \ff{|X_s-Y_s|^2}{\xi_s}\le & \ff{\d M_s}{\xi_s}+ \ff{\kk_2(s)^2\W_2(\mu_s,\nu_s)^2}{2}\d s \\
&-\ff{|X_s-Y_s|^2}{\xi_s^2}\Big\{2-\kk_1(s)\xi_s+ \xi_s'- \ff 1 2\Big\} \d s,\ \ s\in [0, t\land \tau_n].\end{split}\end{equation} By \eqref{Xi} and the monotonicity of $\kk_1$, we have
$$2-\kk_1(s)\xi_s+ \xi_s'- \ff \theta 2\ge 2-\kk_1(T)\xi_s+ \xi_s'- \ff \theta 2=\ff\theta 2.$$
Moreover, since \eqref{A2} implies $(H2)$ for $K_1 = \ff{\kk_1 +\ss{\kk_1^2+ 4\kk_2^2}}2\le \kk_1+\kk_2,$ it follows from Theorem \ref{T1.1} that
$$\W_2(\mu_s,\nu_s)\le W_2(\mu_0,\nu_0)\e^{s\{\kk_1(T)+\kk_2(T)\}},\ \ s\in [0,T].$$ Substituting these into \eqref{V1} and using \eqref{I1}, we arrive at
\beq\label{V2} \beg{split} &\E_{\Q_{t,n}} \int_0^{t\land\tau_n} \ff{|X_s-Y_s|^2}{\xi_s^2}\d s\\
&\le \ff 2 {\xi_0} + \E_{\Q_{t,n}} \int_0^{t\land\tau_n} \kk_2(s)^2 \W_2(\mu_s,\nu_s)^2\d s\\
&\le \bigg[\ff{2}{\xi_0} +  T\kk_2(T)^2\exp[2T(\kk_1(T)+\kk_2(T))] \bigg]\W_2(\mu_0,\nu_0)^2.\end{split}\end{equation}
Writing
$$\log R_{t\land \tau_n}= \int_0^{t\land \tau_n} \ff 1 {\xi_s} \big\<\si_s(X_s)^{-1}(Y_s-X_s),\d \tt W_s\big\>+ \ff 1 2 \int_0^{t\land\tau_n}  \ff {|\si_s(X_s)^{-1}(Y_s-X_s)|^2} {\xi_s^2}\d s,$$ by $\|\si_t^{-1}\|\le \ll(t)$ due to  \eqref{A1}, $\xi_0= \ff{1}{\kk_1(T)} (1-\e^{-\kk_1(T)T})$ due to \eqref{Xi}, and using \eqref{V2}, we arrive at
\beg{align*} \E[R_{t\land \tau_n}\log R_{t\land \tau_n}]&= \ff 1 2 \E_{\Q_{t,n}}\int_0^{t\land\tau_n}  \ff {|\si_s(X_s)^{-1}(Y_s-X_s)|^2} {\xi_s^2}\d s \le \phi (0,T)\W_2(\mu_0,\nu_0)^2.\end{align*} Therefore, \eqref{ES'} holds since $t\in (0,T)$ and $n\ge 1$ are arbitrary.
 \end{proof}

\beg{proof}[Proof of Theorem \ref{T3.1}]  (1) By Lemma \ref{L3.1} and the Girsanov theorem, $\d\Q:=R_T\d\P$ is a probability measure such that
\beq\label{WT} \tt W_s:= W_s - \int_0^s\ff { \si_t(X_t)^{-1}(Y_t-X_t)} {\xi_t} \d t,\ \ s\in [0,T]\end{equation} is a $d$-dimensional Brownian motion.   Then \eqref{CY} reduces to
\beq\label{ET}\d Y_t =  b_t(Y_t,\nu_t)+   \si_t(Y_t) \d \tt W_t.\end{equation}   Consider the DDSDE
$$\d\tt X_t= b_t(\tt X_t, \scr L_{\tt X_t}|_{\tt \P})\d t+ \si_t(\tt X_t)\d \tt W_t,\ \  \tt X_0= Y_0.$$ By the weak uniqueness we  have
$\scr L_{\tt X_t}|_{\tt \P}= P_t^*\nu_0= \nu_t$ for $t\in [0,T]$. Combining this with   \eqref{ET} and the strong uniqueness, we conclude that
$\tt X_t=Y_t$ for $t\in [0,T]$. In particular, $\L_{Y_T}= \nu_T$ as required in (S2). Therefore,  \eqref{LHI} and Lemma \ref{L3.1} lead to
\beg{align*}  (P_T \log f)(\nu_0)\le \log (P_Tf)(\mu_0) + \phi (0,T)\W_2(\mu_0,\nu_0)^2.\end{align*}
In particular,
$$P_T \log f(x) \le (\log P_Tf)(y)+ \phi(0,T) |x-y|^2,\ \ x,y\in \R^d, f\in \B_b(\R^d).$$ According to \cite[Proposition 2.3]{ATW2}, this implies \eqref{LHH0}.

(2) Let $\W_2(\mu_0,\nu_0)>0$. We first assume that $\mu_0$ is absolutely continuous with respect to the Lebesgue measure.   In this case, by \cite[Theorem 10.4.1]{V} (see   \cite{M} when $\nu_0$ is also absolutely continuous),  there exists a measurable map $F: \R^d\to \R^d$ such that $\Xi(x):= x+F(x)$ maps $\mu_0$ into $\nu_0$; that is,
$\nu_0= \mu_0\circ \Xi^{-1}.$ Let
$$\Xi_s(x)=x+sF(x),\ \ \mu_s= \mu_0\circ \Xi_s^{-1},\ \ \ s\in [0,1].$$ Then it is easy to see that
$$\W_2(\mu_s,\mu_t)= |t-s|\W_2(\mu_0,\nu_0),\ \ s,t\in [0,1].$$ Now, for any $n\ge 1$ and $0\le i\le n-1$, we have
$$\W_2(\mu_{i/n}, \mu_{(i+1)/n})=\vv_n:= \ff 1 n \W_2(\mu_0. \nu_0).$$ For any $f\in \B_b(\R^d)$ and $c>0$, when $n$ is large enough such that
$c\vv_n f+1>0$,    the log-Harnack inequality implies
\beq\label{LL}  P_T\log (c\vv_n f+1)(\mu_{i/n})\le \log (c\vv_n P_Tf +1)(\mu_{(i+1)/n}) + \vv_n^2 \phi (0,T),\ \ 0\le i\le n-1.\end{equation}
By Taylor's expansion, there exists a constant $c(f)>0$ depending on $\|f\|_\infty$  such that
\beg{align*} & \Big|P_T\log (c\vv_n f+1)(\mu_{i/n})- c\vv_n (P_T f)(\mu_{i/n})+\ff {(c\vv_n)^2} 2 (P_Tf^2)(\mu_{i/n})\Big|\le \ff {c(f)}{n^3},\\
&\Big|\log (c\vv_n P_Tf +1)(\mu_{(i+1)/n})- c\vv_n (P_T f)(\mu_{(i+1)/n})+\ff {(c\vv_n)^2} 2 (P_Tf)^2(\mu_{(i+1)/n})\Big|\le \ff {c(f)}{n^3}.\end{align*}
Substituting these into \eqref{LL}, we obtain
\beg{align*} &\big|(P_T f)(\mu_{(i+1)/n})- (P_T f)(\mu_{i/n})\big| \\
&\le \ff {c\vv_n} 2\big|(P_Tf^2)(\mu_{i/n}) -(P_T f)^2(\mu_{(i+1)/n})\big|+ \ff{\phi (0,T)\vv_n}c
+\ff{2c(f)}{cn^2\W_2(\mu_0,\nu_0)},\ \ 0\le i\le n-1.\end{align*} Therefore,
\beg{align*} &|(P_T f)(\mu_0)-(P_Tf)(\nu_0)|\le \sum_{i=0}^{n-1}\big|(P_T f)(\mu_{(i+1)/n})- (P_T f)(\mu_{i/n})\big| \\
 &\le \ff{c}2 \sum_{i=0}^{n-1} \vv_n \big|(P_Tf^2)(\mu_{i/n}) -(P_T f)^2(\mu_{(i+1)/n})\big| + \ff{\phi (0,T)\W_2(\mu_0,\nu_0)}c+ \text{O}(n^{-1}).\end{align*} Noting that $\vv_n= \ff 1 n \W_2(\mu_0,\nu_0)$,
by letting $n\to\infty$, we obtain
$$\ff{|(P_T f)(\mu_0)-(P_Tf)(\nu_0)|}{\W_2(\mu_0,\nu_0)} \le \ff{c}2  \sup_{\nu\in B(\mu_0,\W_2(\mu_0,\nu_0))}  \big|(P_Tf^2)(\nu) -(P_T f)^2(\nu)\big| + \ff{\phi (0,T) }c.$$
Minimizing the upper bound in $c>0$, we prove \eqref{LHH1}. Since
$$(P_T^*\nu)(A)- \{(P_T^*\nu)(A)\}^2\le \ff 1 4,\ \  A\in \B(\R^d), \nu\in \scr P,$$  \eqref{LHH2} follows from  \eqref{LHH1} with $f=1_A$.

In general, for any $\mu_0\in \scr P_2$, we take a sequence $\{\mu_0^{(n)}\}_{n\ge 1}\subset \scr P_2$ converging to $\mu_0$ in $\W_2$ and having densities with respect to the Lebesgue measure. Then $\mu_0^{(n)}$ converges to $\mu_0$ weakly. By \eqref{LHH0}, this implies
\beg{align*} &\lim_{n\to 0} (P_{s,t}f(\mu_0^{(n)})=  (P_{s,t}f(\mu_0),\\
 &\lim_{n\to \infty} \sup_{\nu\in B(\mu_0^{(n)}, \W_2(\mu_0^{(n)},\nu_0))}\big\{(P_{s,t} f^2)(\nu)- (P_{s,t} f)^2(\nu)\big\}\\
   &=\sup_{\nu\in B(\mu_0, \W_2(\mu_0,\nu_0))}\big\{(P_{s,t} f^2)(\nu)- (P_{s,t} f)^2(\nu)\big\}.\end{align*} Therefore, by \eqref{LHH1} with $\mu_0^{(n)}$ replacing $\mu_0$ which we just proved, and letting $n\to\infty$, we finish the proof.
\end{proof}

\subsection{Proof of  Theorem \ref{T3.2}}
Again, we only prove for $s=0$. By \eqref{VD}, for any $r>0$  we have
\beq\label{DS}\beg{split} &\exp\bigg[-\ff{r\GG(T)\W_2(\mu_0,\nu_0)^2}{\ll(T)^2}\bigg] \E_{\Q}\Big[ \e^{r\int_0^{t\land\tau_n} \ff{|X_s-Y_s|^2}{\xi_s^2}\d s}\Big] \\
&\le \E_{\Q} \Big[\e^{\ff{2r|X_0-Y_0|^2}{\xi_0} + 4r\int_0^{t\land\tau_n} \ff{\<X_s-Y_s,\{\si_s(X_s)-\si_s(Y_s)\}\d \tt W_s\>}{\xi_s}}\Big]\\
&\le \E_{\Q} \Big[\e^{\ff{2r|X_0-Y_0|^2}{\xi_0}}\E_{\Q}\Big(\e^{4r\int_0^{t\land\tau_n} \ff{\<X_s-Y_s,\{\si_s(X_s)-\si_s(Y_s)\}
\d \tt W_s\>}{\xi_s}}\Big|\F_0\Big)\Big]\\
&\le \E_{\Q} \bigg[\e^{\ff{2r|X_0-Y_0|^2}{\xi_0}}\ss{\E_{\Q}\Big(\e^{32r^2\int_0^{t\land\tau_n} \ff{|\{\si_s(X_s)-\si_s(Y_s)\}^*(X_s-Y_s)|^2}{\xi_s^2}\d s } \Big|\F_0\Big)}\bigg]\\
&\le \ss{ \E_{\Q} \e^{\ff{4r|X_0-Y_0|^2}{\xi_0}}}
\ss{\E_{\Q} \e^{32r^2\gg(T)^2\int_0^{t\land\tau_n} \ff{|X_s-Y_s|^2}{\xi_s^2}\d s}},\end{split}\end{equation}
where   we have used the inequality
$\E_\Q(\e^{M_t}|\F_0)\le \ss{\E_\Q(\e^{2\<M\>_t} |\F_0)}$ for a continuous $\Q$-martingale $M_t$. When $r\le \ff 1 {32 \gg(T)^2}$,
by   Jensen's inequality
$$\ss{\E_{\Q} \e^{32r^2\gg(T)^2\int_0^{t\land\tau_n} \ff{|X_s-Y_s|^2}{\xi_s^2}\d s}}\le \bigg(\E_{\Q}\Big[ \e^{r\int_0^{t\land\tau_n} \ff{|X_s-Y_s|^2}{\xi_s^2}\d s}\Big]\bigg)^{16r\gg(T)^2},$$ so that \eqref{DS} implies
$$ \E_\Q \Big[\e^{r\int_0^{t\land\tau_n} \ff{|X_s-Y_s|^2}{\xi_s^2}\d s}\Big] \le \e^{\ff{r\GG(T)\W_2(\mu_0,\nu_0)^2}{\ll(T)^2(1-16r\gg(T)^2)}} \bigg(\E_\Q \Big[\e^{\ff{4r\kk_1(T)|X_0-Y_0|^2}{1-\e^{-\kk_1(T)T}}}\Big]\bigg)^{\ff 1{2-32r\gg(T)^2}}.$$
Letting $n\uparrow\infty$ and $t\uparrow T$, and noting that $\Q|_{\F_0}=\P|_{\F_0}$ since $R_0=1$, we obtain
\beq\label{ES41}\beg{split}   \E_\Q \Big[\e^{r\int_0^{T} \ff{|X_s-Y_s|^2}{\xi_s^2}\d s}\Big] \le &\e^{\ff{r\GG(T)\W_2(\mu_0,\nu_0)^2}{\ll(T)^2(1-16r\gg(T)^2)}} \bigg(\E \Big[\e^{\ff{4r\kk_1(T)|X_0-Y_0|^2}{1-\e^{-\kk_1(T)T}}}\Big]\bigg)^{\ff 1{2-32r\gg(T)^2}},\\
&\  \ {\rm if}\ 0\le r\le \ff 1{32\gg(T)^2}.\end{split}\end{equation}

On the other hand, letting $M_t= \int_0^t \ff 1 {\xi_s}\<\si_s(X_s)^{-1}(Y_s-X_s), \d\tt W_s\>,\ t\in [0,T],$ by \eqref{R} we have
$R_T= \e^{M_T+\ff 1 2\<M\>_T}$. So,
\beq\label{LT}\beg{split} &\E R_T^{\ff p{p-1}}= \E_\Q R_T^{\ff 1 {p-1}} =\E_\Q \e^{\ff{M_T}{p-1} + \ff{\<M\>_T}{2(p-1)}}\\
&\le \bigg(\E_\Q\exp\bigg[\ff{1}{\ss p-1} M_T-\ff{1}{2(\ss p-1)^2} \<M\>_T\bigg]\bigg)^{\ff 1 {1+\ss p}}
\bigg(\E_\Q\exp\bigg[\ff{\<M\>_T}{2(\ss p-1)^2}\bigg]\bigg)^{\ff{\ss p}{\ss p+1}}\\
&\le \bigg(\E_\Q\exp\bigg[\ff{\ll(T)^2}{2(\ss p-1)^2}\int_0^T \ff{|X_s-Y_s|^2}{\xi_s^2}\d s\bigg]  \bigg)^{\ff{\ss p}{\ss p+1}}.\end{split}\end{equation}
Since $p\ge p(T)$ implies
$\ff{\ll(T)^2}{2(\ss p-1)^2}\le \ff 1 {32 \gg(T)^2},$ this and \eqref{ES41} with $r= \ff{\ll(T)^2}{2(\ss p-1)^2}$  yield
\beg{align*}  \E R_T^{\ff p{p-1}}
 \le \e^{\ff{\ss p\,\GG(T) \W_2(\mu_0,\nu_0)^2}{(\ss p+1)\{2(\ss p-1)^2-16\ll(T)^2\gg(T)^2\}}}
\bigg(\E_\Q\e^{\ff{2\ll(T)^2\kk_1(T)|X_0-Y_0|^2}{(\ss p-1)^2 (1-\e^{-\kk_1(T)T})}}\bigg)^{\ff{\ss p\,(\ss p-1)^2 }{(\ss p+1)\{2(\ss p-1)^2-16\ll(T)^2\gg(T)^2\}}}.\end{align*}
Substituting into \eqref{HI}, we finish the proof.

\section{Shift Harnack inequality and applications}

In this section we establish the shift Harnack inequality and integration by parts formula introduced in \cite{W14a}. Since the study for the multiplicative noise case is very complicated, here we only consider the additive noise for which the DDSDE \eqref{E1} reduces to
\beq\label{E5} \d X_t= b_t(X_t, \L_{X_t})\d t +\si_t \d W_t.\end{equation}

\beg{thm}\label{T5.1} Let $\si: [0,\infty)\to \R^d\otimes \R^d$ and $b: [0,\infty)\times \R^d\times\scr P_2\to \R^d$ are measurable such that
$\si_t$ is invertible with $\|\si_t\|+\|\si_t^{-1}\|$  locally bounded in $t\ge 0$, and $b_t(\cdot,\mu_t)$ is differentiable with
$$\int_0^T \|\nn b_t(\cdot, \mu_t)\|_\infty^2\d t<\infty,\ \ T>0, \mu_\cdot\in C([0,T]; \scr P_2).$$
\beg{enumerate} \item[$(1)$] For any $p>1, t>s\ge 0, \mu_0\in \scr P_2, v\in\R^d$ and $f\in \B_b^+(\R^d)$,
\beg{align*}(P_{s,t}f)^p(\mu_0)\le &(P_{s,t}f^p(v+\cdot))(\mu_0)\\
&\times \exp\bigg[\ff{p\ss p\, |v|^2\int_s^t \|\si_r^{-1}\|^2 \big\{1
+ (r-s)\|\nn b_r(\cdot,P_{s,r}^*\mu_0)\|_\infty \big\}^2\d r}{2(p-1)(\ss p+1)(t-s)^2}\bigg].\end{align*} Moreover, for any $f\in \B_b^+(\R^d)$,
$$(P_{s,t}\log f)(\mu_0)\le \log (P_{s,t} f(v+\cdot))(\mu_0)+ \ff {|v|^2}{2(t-s)^2} \int_s^t \|\si_r^{-1}\|^2\Big( 1+ (r-s)\|\nn b_r(\cdot, \mu_r)\|_\infty\Big)^2\d r.$$
\item[$(2)$] For any $t>s\ge 0, f\in C^1(\R^d)$ and $\F_s$-measurable random variable $X_{s,s}$ with $\mu_0:=\L_{X_{s,s}}\in \scr P_2$,
$$\E(\nn_v f)(X_{s,t})= \E\bigg[\ff{f(X_{s,t})}{t-s}\int_s^t (r-s)\big\<\si_r^{-1} \nn_v b_r(\cdot, P_{s,r}^*\mu_0)(X_{s,r}),\ \d W_r\big\>\bigg],\ \ v\in \R^d.$$
\end{enumerate} \end{thm}

\beg{proof} Without loss of generality, we only prove for $s=0$ and $t=T$ for some fixed time $T>0$.
Denote $\mu_t=P_t^*\mu_0=\L_{X_t}, t\ge 0$. Then \eqref{E5} becomes
\beq\label{E5'}  \d X_t= b_t(X_t, \mu_t)\d t +\si_t \d W_t,\ \ \L_{X_0}=\mu_0.\end{equation}
Let $Y_t=X_t+\ff {tv}T,\ t\in [0,T]$. Then
$$ \d Y_t= b_t(Y_t, \mu_t)\d t +\si_t \d \tt W_t,\ \ \L_{Y_0}=\mu_0, t\in [0,T], $$ where
\beg{align*} &\tt W_t:= W_t +\int_0^t \xi_s\d s,\\
&\xi_t:= \si_t^{-1}\Big\{\ff v T + b_t(X_t,\mu_t)-b_t\Big(X_t+\ff {tv}T, \mu_t\Big)\Big\}.\end{align*}
Let $R_T= \exp[-\int_0^T \<\xi_t, \d W_t\>-\ff 1 2\int_0^T |\xi_s|^2\d s].$ By the Girsanov theorem we obtain
$$(P_T f)(\mu_0)=\E [R_T f(Y_T)]= \E[R_T f(X_T+v)]\le (P_T f^p(v+\cdot))^{\ff 1 p}(\mu_0) \big(\E R_T^{\ff p{p-1}}\big)^{\ff {p-1}p}.$$
 This proves  (1) since similarly to \eqref{LT}, we have
\beg{align*} &\E R_T^{\ff p {p-1}}= \E_\Q R_T^{\ff 1 {p-1}}
 \le \Big(\E_\Q\e^{\ff{p}{2(p-1)^2} \int_0^T |\xi_s|^2\d s}\Big)^{\ff{\ss p}{\ss p+1}} \\
&\le \exp\bigg[\ff{p\ss p\, |v|^2\int_0^T \|\si_t^{-1}\|^2 \big\{1
+ t\|\nn b_t(\cdot,P_{t}^*\mu_0)\|_\infty \big\}^2\d t}{2(p-1)^2(\ss p+1)T^2}\bigg].\end{align*}

To prove (2), we let $X_t^\vv= X_t+\ff{\vv tv}{T}$ for $\vv\in (0,1)$ and $t\in [0,T]$. Using $\vv v$ replace $v$, the above argument implies
$$(P_T f)(\mu_0)= \E[R_T^\vv f(X_T+ \vv v)],\ \ \vv\in (0,1),$$
where
\beg{align*} & R_T^\vv:= \exp\bigg[-\int_0^T \<\xi_t^\vv, \d W_t\>-\ff 1 2\int_0^T |\xi_s^\vv|^2\d s\bigg],\\
& \xi_t^\vv:=  \si_t^{-1}\Big\{\ff {\vv v} T + b_t(X_t,\mu_t)-b_t\Big(X_t+\ff {\vv tv}T, \mu_t\Big)\Big\}.\end{align*}
Therefore,
\beg{align*} 0&= \lim_{\vv\to 0} \ff 1 \vv \E[R_T^\vv f(X_T+ \vv v)- f(X_T)]\\
&= \E[(\nn_vf)(X_T)] - \E\bigg[\ff{f(X_{T})}{T}\int_0^T  r \big\<\si_r^{-1} \nn_v b_r(\cdot, P_{0,r}^*\mu_0)(X_{r}),\ \d W_r\big\>\bigg].\end{align*}
Then the proof is finished.
\end{proof}

As applications of Theorem \ref{T5.1}, we have the following estimates on the  density of $P_{s,t}^*$.

\beg{cor} \label{C5.2} In the situation of Theorem $\ref{T5.1}$, for any $t>s\ge 0$ and $\mu_0\in \scr P_2$, $(P_{s,t}^*\mu_0)(\d x)= \rr_{s,t}^{\mu_0}(x)\d x$
for some   density function $\rr_{s,t}^{\mu_0}$ satisfying the following estimates:
\beq\label{ET1} \beg{split}&\int_{\R^d} |\nn \log \rr_{s,t}^{\mu_0}(x)|^p\rr_{s,t}^{\mu_0}(x)\d x\\
&\le  \bigg\{\ff{ (1\lor\ff{p(p-1)}2 )  }{(t-s)^2} \int_s^t (r-s)^2 \|\si_r^{-1}\|^2 \|\nn b_r(\cdot, P_{s,r}^*\mu_0)\|_\infty^2 \d r\bigg\}^{\ff p 2 (1\land \ff 1 {p-1})},\ \ p>1;\end{split}\end{equation}
 \beq\label{ET2} \beg{split} &\int_{\R^d} \big\{\rr_{s,t}^{\mu_0}(x)\big\}^{\ff p{p-1}}\d x \\
 &\le \bigg(\ff{p\ss p\int_s^t \|\si_r^{-1}\|^2\{1+(r-s)\|\nn b_2(\cdot, P_{s,r}^*\mu_0)\|_\infty\}^2\d r}{4\pi (p-1)(\ss p+1)(t-s)^2}\bigg)^{\ff d{2(p-1)}},\ \ p>1;\end{split}\end{equation}
\beq\label{ET3}  \beg{split} &\int_{\R^d} \rr_{s,t}^{\mu_0}(x)\log  \rr_{s,t}^{\mu_0}(x)\d x\\
&\le \ff d 2 \log \bigg(\ff{\int_s^t \|\si_r^{-1}\|^2\{1+(r-s)\|\nn b_2(\cdot, P_{s,r}^*\mu_0)\|_\infty\}^2\d r}{4\pi(\ss p+1)(t-s)^2}\bigg).\end{split}\end{equation}
\end{cor}

\beg{proof} According to \cite[Theorem 2.4, Theorem 2.5]{W14a}, the desired assertions follow from Theorem \ref{T5.1} . Indeed, by \cite[Theorem 2.4]{W14a}, the integration by parts formula in Theorem \ref{T5.1}(2) implies the existence of density $\rr_{s,t}^{\mu_0}$ and
$$\int_{\R^d} |\nn \log \rr_{s,t}^{\mu_0}(x)|^p\rr_{s,t}^{\mu_0}(x)\d x = \E |\nn \log \rr_{s,t}^{\mu_0}|^p(X_{s,t})
 = \E\big|\E(N|X_{s,t})\big|^p\le \E|N|^p,$$ where
 $$N:= \ff{1}{t-s}\int_s^t (r-s) \{\si_r^{-1} \nn  b_r(\cdot, P_{s,r}^*\mu_0)(X_{s,r})\}^*  \d W_r.$$
Then estimate \eqref{ET1} follows.  Moreover, it is easy to see that the proof of \cite[Theorem 2.5]{W14a} also applies to $P_{s,t}^*\mu_0$ in place of $P(x,\cdot)$, so that estimates \eqref{ET2} and \eqref{ET3} follow from Theorem \ref{T5.1}(1).

\end{proof}

\section{DDSDEs for homogeneous Landau equation}

We consider  the homogeneous Landau equation with $r\in [0,1]$    on $\R^3$ (see e.g. \cite{V2}):
\beq\label{LLM}\pp_t f_t= \ff 1 2 {\rm div}\bigg(\int_{\R^3} |\cdot-y|^{2+\gg} \Big(I- \ff{(\cdot-y)\otimes (\cdot-y)}{|\cdot-y|^2}\Big)\Big\{f_t(y)\nn f_t-f_t\nn f_t(y)\Big\}\d y\bigg).\end{equation} Let $a(x)= |x|^\gg(|x|^2I-x\otimes x)$ and
\beq\label{LM} b_0(x):={\rm div} a(x)= -2|x|^\gg x,\ \
 \si_0(x):= |x|^{\ff \gg 2}\left(\beg{matrix} &x_2 &0 &x_3\\
&-x_1 &x_3 &0\\
&0 &-x_2 &-x_1\end{matrix}\right).\end{equation} Then $\si_0\si_0^*=a$. Take
\beq\label{BS} \beg{split} &b_t(x,\mu)= b(x,\mu):=-2\int_{\R^d} |x-z|^\gg (x-  z)\mu(\d z),\\
&\si_t(x,\mu)= \si(x,\mu):=\int_{\R^d}   \si_0(x-z)\mu(\d z).\end{split}\end{equation}
Then the density of $\scr L_{X_t}$ for the DDSDE \eqref{E1} is a weak solution to \eqref{LLM}.
In this section we consider \eqref{E1} for this specific choice of $b$ and $\si$.

\subsection{The case with Maxwell molecules: $\gg=0$ }
When $\gg=0$, both  $b_0$ and $\si_0$ in \eqref{LM} are Lipschitz continuous. Below we consider  a more general model. For two Lipschitz continuous maps
$$b_0: \R^d \to \R^d,\ \ \si_0: \R^d\to \R^d\otimes \R^d,$$  let
$$ b^\aa(x,\mu):= \int_{\R^d} b_0(x-\aa z)\mu(\d z),\ \ \si^\aa(x,\mu):= \int_{\R^d} \si_0(x-\aa z)\mu(\d z),\ \ \ \ \aa\in\R, x\in\R^d, \mu\in \scr P_2.$$ For  fixed $\aa,\bb\in \R$,  consider the DDSDE
\beq\label{EP2} \d X_t = b^\aa(X_t, \L_{X_t})\d t + \si^\bb(X_t, \L_{X_t})\d W_t.\end{equation}

\beg{thm}\label{TN} Let $\aa,\bb\in \R$, $B_0:=\|\nn b_0\|_\infty<\infty$ and $C_0:=\sup_{|v|=1,x\in\R^d} \|(\nn_v\si_0)(x)\|_{HS}^2<\infty,$ where $\nn_v$ denotes the directional derivative along $v$. Moreover, assume that 
$$\<b_0(x)-b_0(y), x-y\>\le K_0|x-y|^2,\ \ x,y\in \R^d $$ hold for some constant $K_0$, which can be negative. 
\beg{enumerate} \item[$(1)$] For any $\F_0$-measurable $X_0$ with $\E |X_0|^2<\infty$,  the equation $\eqref{EP2}$ has a unique solution $(X_t)_{t\ge 0}$, and $\sup_{t\in [0,T]}\E |X_t|^2<\infty$ for all $T>0$. Moreover, $X_t(x)$ is jointly continuous in $(t,x)\in [0,\infty)\times \R^d$, where $(X_t(x))_{t\ge 0}$ is the unique solution with $X_0=x$. 
\item[$(2)$]  For any $\mu_0,\nu_0\in \scr P_2$,
$$\W_2(P_t^*\mu_0, P_t^*\nu_0)^2\le \W_2(\mu_0, \nu_0)^2\e^{(2K_0 + C_0(1+|\bb|)^2+2|\aa|B_0)t},\ \ t\ge 0.$$
If, in particular, $2K_0 + C_0(1+|\bb|)^2+2|\aa|<0$, then $P_t^*$ has a unique invariant probability measure.
\item[$(3)$] If $\bb=0$ and $\si_0$ is invertible with $\ll:= \|\si_0^{-1}\|_\infty<\infty$, then assertions in Theorem $\ref{T3.1}$, Theorem $\ref{T3.2}$ and Corollary $\ref{C3.3}$ hold for
    $$\phi(s,t)= \ll^2\bigg(\ff{2K_0+C_0}{1-\e^{-(2K_0+C_0)(t-s)}}+  (t-s)|\aa|\e^{2(t-s)(2K_0+C_0+2|\aa|)} \bigg).$$
\item[$(4)$] If $\si_0$ is constant and invertible, then assertions in Theorem $\ref{T5.1}$ and Corollary $\ref{C5.2}$ hold  for $\si_r\equiv \si_0$.
\end{enumerate} \end{thm}

\beg{proof} Since $b_0$ and $\si_0$ are Lipschitz continuous, it is easy to see that $(H1)$-$(H3)$ and \eqref{X1}  hold for $(b_t,\si_t)\equiv (b^\aa, \si^\bb)$ for all $t\ge 0$. Then the first assertion follows from Theorems \ref{T1.1} and \ref{T2.2}.

To prove the second assertion using Theorem \ref{T1.2}(2), we observe that for any $\pi\in \C(\mu,\nu)$.
\beg{align*} &\<b^\aa(x,\mu)-b^\aa(y,\nu), x-y\>\\
&= \int_{\R^d}\Big(\<b_0(x-\aa z)-b_0(y-\aa z), x-y\>+ \<b_0(y-\aa z)-b_0(y-\aa z'),x-y\>\Big)\pi(\d z, \d z')\\
&\le K_0|x-y|^2 +B_0|\aa|\cdot |x-y|\int_{\R^d} |z-z'|\pi(\d z, \d z').\end{align*} Then
\beq\label{BS2}\beg{split}  2 \<b^\aa(x,\mu)-b^\aa(y,\nu), x-y\>&\le 2K_0|x-y|^2 + 2 |\aa| B_0 \W_1(\mu,\nu)|x-y|\\
&\le (2K_0+|\aa|B_0)|x-y|^2 + |\aa|B_0 \W_2(\mu,\nu)^2.\end{split}\end{equation}
Similarly,
\beg{align*}  \|\si^\bb(x,\mu)-\si^\bb(y,\nu)\|_{HS}^2&\le C_0 \big\{|x-y|+ |\bb| \W_1(\mu,\nu)\big\}^2\\
&\le C_0(1+|\bb|)|x-y|^2+ C_0(|\bb|+\bb^2)\W_2(\mu,\nu)^2.\end{align*}
Combining this with \eqref{BS2} we obtain
\beg{align*} &2 \<b^\aa(x,\mu)-b^\aa(y,\nu), x-y\>+\|\si^\bb(x,\mu)-\si^\bb(y,\nu)\|_{HS}^2\\
&\le \{2K_0+|\aa|B_0+C_0(1+|\bb|)\}|x-y|^2+\{|\aa|B_0+C_0|\bb|(1+|\bb|)\}\W_2(\mu,\nu)^2.\end{align*}
Then the second assertion follows from Theorem \ref{T1.2}(2).

Finally, by \eqref{BS2} and $\|\si_0(x)-\si_0(y)\|_{HS}^2\le C_0|x-y|^2$ we have
$$2 \<b^\aa(x,\mu)-b^\aa(y,\nu), x-y\>+\|\si_0(x)-\si_0(y)\|_{HS}^2\le (2K_0+ C_0 )|x-y|^2+ 2|\aa| |x-y|\W_2(\mu,\nu).$$ Then assumption {\bf (A)} holds for
$\ll(t)=\ll, \kk_1(t)= 2K_0+C_0$ and $\kk_2(t)= 2|\aa|.$ Therefore, assertions (3) and (4) follow from Theorem $\ref{T3.1}$, Theorem $\ref{T3.2}$, Corollary $\ref{C3.3}$,  Theorem $\ref{T5.1}$ and Corollary $\ref{C5.2}$.
\end{proof}

Coming back to the DDSDE for the homogeneous Landau equation with Maxwell molecules, i.e. $\eqref{EP2}$ for  $b_0$ and $\si_0$ in \eqref{LM}, Theorem \ref{TN} applies with $B_0=C_0=2$ and $K_0=-2$, so that we have  the following result.

\beg{cor}\label{CC} Let $b_0$ and $\si_0$ be in $\eqref{LM}$ and let $\aa,\bb\in \R$. For any $\F_0$-measurable $X_0$ with $\E |X_0|^2<\infty$,  the equation $\eqref{EP2}$ has a unique solution and $\sup_{t\in [0,T]}\E |X_t|^2<\infty$ for all $T>0$. Moreover, $X_t(x)$ is jointly continuous in $(t,x)\in [0,\infty)\times \R^d$.
 Moreover, for any $\mu_0,\nu_0\in \scr P_2$,
$$\W_2(P_t^*\mu_0, P_t^*\nu_0)^2\le \W_2(\mu_0, \nu_0)^2\e^{\{4(|\aa|+|\bb|)+2\bb^2-2\}t},\ \ t\ge 0,\mu_0,\nu_0\in \scr P_2.$$
When $2(|\aa|+|\bb|)+\bb^2<1$, $P_t^*$ has a unique invariant probability measure $\mu$ and
$$\W_2(P_t^*\nu_0, \mu)^2\le \e^{-2(1-2|\aa|-2|\bb|-\bb^2)t}\W_2(\nu_0,\mu)^2,\ \ t\ge 0, \nu_0\in \scr P_2.$$
When $\bb=\aa=1$ which corresponds to the homogeneous Landau equation with Maxwell molecules,
\beq\label{*CC0}\W_2(P_t^*\mu_0, P_t^*\nu_0)^2\le \W_2(\mu_0, \nu_0)^2\e^{8t},\ \ t\ge 0,\mu_0,\nu_0\in \scr P_2.\end{equation}
 \end{cor}

 \paragraph{Remark 6.1.} Let $N(z, A)$ denote the normal distribution on $\R^d$ with mean $z\in\R^d$ and covariance $A$, and let   $\aa=\bb=1$ in Corollary \ref{CC}  for the homogeneous Landau equation with Maxwell molecules.  According to \cite[Theorem 1.1]{Carlen} (see also \cite{DV2}),  there exists a constant $p>0$ such that if
 $$\int_{\R^d} |x|^p \mu_0(\d x)+ \int_{\R^d}|\xi|^p|\hat\mu_0(\xi)|^2\d\xi  <\infty,$$ where $\hat \mu_0$ is the fourier transform of $\mu_0$,  then
 $$\|  P_t^*\mu_0 -N(z_0, \gg_0^2 I)\|_{var}  \le c_1\e^{-c_2t},\ \ t\ge 0$$ holds for some constants $c_1,c_2>0$ depending on $\mu_0$, where $z_0:=\int_{\R^3}x\mu_0(\d x)$ and    $\gg_0^2:=\int_{\R^d} |x-z_0|^2\mu_0(\d x).$ See   \cite{CA} for exponential convergence in the case that $\gg\in (0,1].$
Therefore, $P_t^*$ is not ergodic since the limit distribution varies in the initial one. This fits the inequality \eqref{*CC0} where the upper bound  does not go to 0 as $t\to\infty$. However, it seems that the sharp upper bound in \eqref{*CC0} should be bounded in $t$.

\subsection{ The case with hard potentials: $\gg\in [0,1]$ }
When $\gg\in [0,1]$,  the weak existence and uniqueness have been proved in \cite{FG}. To prove the same assertion  for strong solutions,  we first present a result for the equivalence of the weak  existence/uniqueness and the strong existence/uniqueness.

\beg{thm} \label{T7.1} Let $\theta\ge 1.$ Assume that for any $\mu\in C([0,\infty)\to \scr P_\theta)$ the SDE
\beq\label{ESDE} \d X_t= b_t(X_t,\mu_t)\d t +\si_t(X_t,\mu_t)\d W_t\end{equation} has strong existence and uniqueness for $X_0$ with $\L_{X_0}=\mu_0$. Then for initial distribution  $\mu_0\in \scr P_\theta,$ the DDSDE $\eqref{E1}$ has weak existence (respectively uniqueness) if and only if it has strong existence $($respectively uniqueness$)$. \end{thm}

\beg{proof} (a) Since  the strong existence implies the weak one, it suffices to prove the strong existence from the weak one.
For any initial distribution $\mu_0\in \scr P_\theta$, let $(\bar X_t, \bar W_t)$ be a weak solution under probability $\bar \P$.  We have
\beq\label{XY1} \d \bar X_t= b_t(\bar X_t, \mu_t) \d t +\si_t(\bar X_t, \mu_t)\d \bar W_t,\end{equation} where $\mu_t:= \L_{\bar X_t}|_{\bar \P}.$ Now, given a
Brownian motion under the probability $\P$, let $X_t$ be a strong solution to \eqref{ESDE} with $\L_{X_0}=\mu_0$. By Yamada-Watanabe's principle for SDE, the strong existence and uniqueness of \eqref{ESDE} imply the weak uniqueness, so that $\L_{X_t}=\mu_t$ so that \eqref{ESDE} reduces to the DDSDE \eqref{E1}. Then  the strong solution to   \eqref{ESDE} is also a strong solution to \eqref{E1}.

(b) Obviously,  the weak uniqueness implies the strong uniqueness. On the other hand, let \eqref{E1} has strong uniqueness, we aim to prove the weak uniqueness. Let $(X_t^{(i)}, W_t^{(i)})$ under probability $\P^i (i=1,2)$  be two weak solutions to \eqref{E1} with $\L_{X_0^{(1)}}|_{\P^1}= \L_{X_0^{(2)}}|_{\P^2}=\mu_0$, we aim to prove
\beq\label{WU} \L_{X^{(1)}}|_{\P^1}=  \L_{X^{(2)}}|_{\P^2}.\end{equation}
Let $\mu_t= \L_{X_t^{(1)}}|_{\P^1}.$ By assumption, the SDE
\beq\label{XY0} \d X_t= b_t(X_t,\mu_t)\d t +\si_t(X_t,\mu_t)\d W_t^{(2)},\ \ X_0= X_0^{(2)}\end{equation}   has a unique strong solution $X:=(X_t)_{t\ge 0}$.
 By  Yamada-Watanabe's principle, \eqref{XY0} also has weak uniqueness. So,
 \beq\label{XY2}\L_{X}|_{\P^2}= \L_{X^{(1)}}|_{\P^1}.\end{equation} In particular, $ \L_{X_t}|_{\P^2}=\mu_t$, so that $X_t$ is also a strong solution to  \eqref{E1} with   the given Brownian motion $W_t^{(2)}$ replacing $W_t$. Since $X_t^{(2)}$  solves the  same DDSDE, by the strong uniqueness of \eqref{E1} we have $X=X^{(2)}$. Combining this with \eqref{XY2}, we prove \eqref{WU}.
\end{proof}

Now, we consider  the DDSDE \eqref{E1} with $b_t$ and $\si_t$ in \eqref{BS} for $\gg \in (0,1].$

\beg{cor} Let   $b_t$ and $\si_t$ in \eqref{BS} for $\gg \in (0,1].$  Then for any  $\F_0$-measurable $X_0$ with density $f_0$ satisfying \eqref{DST},    the DDSDE $\eqref{E1}$ has a unique strong solution such that
$\E  \e^{|X_t|^\aa} <\infty$ for any $t>0.$\end{cor}

\beg{proof} By \cite[Theorem 2]{FG}, the SDDE has a unique weak solution such that $\E  \e^{|X_t|^\aa} <\infty$ for any $t>0.$ According to Theorem \ref{T7.1}, the same holds for the strong solution. \end{proof}

  \paragraph{Acknowledgement.} The author would like to thank the referees as well as Jianhai Bao, Arnaud Guillin and Xing Huang  for corrections and helpful comments.

\end{document}